\documentclass[11pt, lettersize, notitlepage, leqno, footskip=0.6cm]{article} 
\usepackage{amssymb}
\usepackage{amsthm}
\usepackage{amsmath}	
\usepackage{graphicx}
\usepackage{amscd}
\usepackage{tikz-cd}	
	
\usepackage{thmtools}
\usepackage[none]{hyphenat}	

\usepackage{verbatim}

\newcommand{\bc}{\mathbb C}
\newcommand{\bF}{\mathbb F}

\newcommand{\bz}{\mathbb Z}

\newcommand{\bq}{\mathbb Q}

\newcommand{\bFp}{\mathbb{F}_p}
\newcommand{\bFP}{\ov{\mathbb{F}}_p}

\newcommand{\bFq}{\mathbb{F}_q}
\newcommand{\bFpk}{\mathbb{F}_{p^k}}

\newcommand{\Gal}{\mathrm{Gal}}

\newcommand{\W}{\mathcal W}

\newcommand{\qbar}{\overline {\mathbb{Q}}}

\newcommand{\absq}{\mathrm{Gal}_{\bq}}

\newcommand{\Fr}{\mathrm{Fr}}

\newcommand{\la}{\langle}
\newcommand{\ra}{\rangle}

\newcommand{\lra}{\longrightarrow}

\newcommand{\hra}{\hookrightarrow}

\newcommand{\wti}{\widetilde}
\newcommand{\mf}{\mathfrak}
\newcommand{\mc}{\mathcal}
\newcommand{\mb}{\mathbb}
 
\newcommand{\mr}{\mathrm}
\newcommand{\mfp}{\mathfrak{p}}

\newcommand{\mfl}{\mathfrak{l}}

\newcommand{\al}{\alpha}
\newcommand{\be}{\beta}

\newcommand{\lamb}{\lambda}

\newcommand{\Lamb}{\Lambda}

\newcommand{\what}{\widehat}

\newcommand{\ov}{\overline}

\newcommand{\sub}{\subseteq}

\newcommand{\GalWB}{\Gal(\bq(\W_{B_{\pi}})/\bq)}

\DeclareMathOperator{\Cent}{Cent}

\DeclareMathOperator{\End}{End}
\DeclareMathOperator{\GL}{GL}

\DeclareMathOperator{\Hom}{Hom}

\DeclareMathOperator{\Mat}{Mat}

\DeclareMathOperator{\inv}{inv}

\DeclareMathOperator{\Nm}{Nm}
\DeclareMathOperator{\Nr}{Nrd}
\DeclareMathOperator{\charac}{char}

\DeclareMathOperator{\conn}{conn}
\DeclareMathOperator{\IV}{IV}
\DeclareMathOperator{\Cl}{Cl}

\DeclareMathOperator{\Iso}{Iso}
\DeclareMathOperator{\isog}{isog}
\DeclareMathOperator{\Jac}{Jac}

\newtheorem{Thm}{Theorem}[section]
\newtheorem{Prop}[Thm]{Proposition}
\newtheorem{Lem}[Thm]{Lemma}
\newtheorem{Corr}[Thm]{Corollary}

\newtheorem{Example}[Thm]{Example}

\newtheorem{Def}{Definition}[section]

\declaretheoremstyle[%
  spaceabove=-2pt,%
  spacebelow=8pt,%
  headfont=\normalfont\itshape,%
  postheadspace=1em,%
  qed=\qedsymbol%
]{mystyle} 
\declaretheorem[name={Proof},style=mystyle,unnumbered,
]{prf}

\numberwithin{equation}{section}

\title{Isogenies of certain abelian varieties over finite fields\\ with $p$-ranks zero} 
\author{Steve Thakur}  
\date{\vspace{-3ex}}

\begin{document} 
\maketitle

\begin{abstract}

 \noindent We study the isogenies of certain abelian varieties over finite fields with non-commutative endomorphism algebras with a view to potential use in isogeny-based cryptography. In particular, we show that any two such abelian varieties with endomorphism rings maximal orders in the endomorphism algebra are linked by a cyclic isogeny of prime degree.\end{abstract}

\section{\fontsize{11}{11}\selectfont Introduction}

In this article, we study the isogenies of a certain class of absolutely simple abelian varieties over finite fields with $p$-rank zero and non-commutative endomorphism rings. 

In recent years, abelian varieties with non-commutative endomorphism rings have seen a growing interest for applications to isogeny-based cryptography. This is primarily due to the quantum attack with sub-exponential complexity on isogeny-based cryptosystems of abelian varies with commutative endomorphism rings  ([CJS10]). The cryptosystems based on isogenies of supersingular elliptic curves have gained a lot of attention in the last few years since they have substantially smaller key sizes than other quantum-resistant schemes. However, such abelian varieties are relatively unexplored for dimensions larger than two. To this end, one of the goals in this article is to explore the behavior of the isogenies of abelian varieties of type $\IV(1,d)$ (see Definition 2.1) over finite fields. This is a class of absolutely simple abelian varieties whose endomorphism algebras have centers that are imaginary quadratic fields.

Such abelian varieties exist for every dimension $\geq 3$. For dimensions larger than three, most principally polarized abelian varieties are not Jacobians and hence, are much less promising for cryptographic purposes. But we explore the isogenies of these abelian varieties for all dimensions since they seem interesting in their own right.

\subsection{\fontsize{11}{11}\selectfont Notations and background}

For an abelian variety $A$ over a field $F$, $\End(A)$ denotes its endomorphism ring and $\End^0(A)$ the endomorphism algebra of $A$ over the algebraic closure of the field of definition. By Honda-Tate theory, we have the well-known bijection $$\{\text{Simple abelian varieties over } \bFq \text{ up to isogeny}\}\longleftrightarrow \{\text{Weil } q\text{-integers up to } \absq\text{-conjugacy} \}$$ induced by the map sending an abelian variety to its Frobenius. For a Weil number $\pi$, we write $B_{\pi}$ for the corresponding simple abelian variety over $\bFq$. The dimension of $B_{\pi}$ is given by $$2\dim B = [\bq(\pi):\bq][\End^0(B_{\pi}):\bq(\pi)]^{1/2}.$$ Note that $\End^0(B_{\pi})$ is a central division algebra over $\bq(\pi)$ and hence, $[\End^0(B_{\pi}):\bq(\pi)]^{1/2}$ is an	integer. The characteristic polynomial of $B_{\pi}$ on the Tate representation $V_l(A):= T_l(A)\otimes_{\bz_l}\bq_l$ (for any prime $l\neq p$) is independent of $l$ and is given by \vspace{-0.2cm}$$P_{B_{\pi}}(X):= \prod\limits_{\sigma\in \absq} (X-\sigma(\pi))^{m_{\pi}}$$ where $m_{\pi} = [\End^0(B_{\pi}):\bq(\pi)]^{1/2}.$ We denote by $\W_{B_{\pi}}$ the set of Galois conjugates of $\pi$. Hence, $\bq(\W_{B_{\pi}})$ is the splitting field of $P_{B_{\pi}}(X)$. The Galois group $\Gal(\bq(\W_{B_{\pi}})/\bq)$ is a subgroup of the wreath product $(\bz/2\bz)^g\rtimes S_g$, the Galois group of the generic CM field of degree $2g$, where $g = [\bq(\pi+\ov{\pi}):\bq]$.

\begin{Def} \normalfont An abelian variety $A$ over a field $F$ is \textit{simple} if it does not contain a strict non-zero abelian subvariety. We say $A$ is \textit{absolutely} or \textit{geometrically} simple if the base change $A\times_F \ov{F}$ to the algebraic closure is simple.\end{Def}

\begin{Def} \normalfont An abelian variety $A$ is \textit{iso-simple} if it has a unique simple abelian subvariety up to isogeny.\end{Def}

\noindent We state a few well-known facts about abelian varieties over finite fields. We refer the reader to notes [Oo95] for proofs and further details.

\begin{Prop} For any simple abelian variety $B$ over a finite field $\bFq$, the abelian variety $B\times_{\bFq} \ov{\bF}_q$ is iso-simple. 
\end{Prop}
\vspace{-0.15cm}
With this setup, let $\pi$ be a Weil number corresponding to $B$ and let $\wti{B}$ be the unique simple component (up to isogeny) of the base change $B\times_{\bFq} \ov{\bF}_q$ to the algebraic closure. Let $N$ be the smallest integer such that $\wti{B}$ has a model over the field $\bF _{q^N}$. Then $\wti{B}$ corresponds to the Weil number $\pi^N$ and we have an isogeny \vspace{-0.15cm}$$B\times_{\bFq} \ov{\bF}_q =_{\isog} \wti{B}^{(\dim B)/N}.$$

\begin{Prop} Let $\pi$ be a Weil $q$-integer and let $B_{\pi}$ be the corresponding simple abelian variety over $\bFq$. The following are equivalent:\\
\noindent $1.$ $B_{\pi}\times_{\bFq} \bF _{q^N}$ is simple.\\
$2.$ $\bq(\pi^N) = \bq(\pi)$.\end{Prop}

\vspace{-0.15cm}
\noindent Thus, $B_{\pi}$ is absolutely simple if and only if $\bq(\pi^N) = \bq(\pi)$ for every integer $N$. 

\begin{Prop} Let $\pi$ be a Weil $q$-integer and write $D_{\pi}:= \End^0_{\bFq}(B_{\pi})$. Then $D_{\pi}$ is a central division algebra over $\bq(\pi)$ and its Hasse invariants are given by $$\mr{inv}_{v}(D_{\pi}) = \begin{cases} 0 & \text{ if } v\nmid q.\\
\frac{1}{2} & \text{ if } v \text{ is real.}\\
[\bq(\pi)_{v}:\bq_p]\frac{v(\pi)}{v(q)} & \text{ if } v|q.\end{cases}$$\end{Prop}

\noindent In particular, $\End^0(B_{\pi})$ is commutative if and only if the local degrees $[\bq(\pi)_{v}:\bq_p]$ annihilate the Newton slopes $\frac{v(\pi)}{v(q)}$. For instance, if $B_{\pi}$ is ordinary, the slopes $\frac{v(\pi)}{v(q)}$ are either $0$ or $1$ and hence, $\End^0(B_{\pi})$ is commutative.



\begin{Prop} Let $B_{\pi}$ be a simple abelian variety over a finite field $\bFq$ corresponding to a Weil number $\pi$ and let $l$ be a prime that does not divide $q$. The order $\big|B_{\pi}(\bFq)\big|$ of the group of $\bFq$-points is given by \vspace{-0.1cm}$$\big|B_{\pi}(\bFq)\big| = P_{B_{\pi}}(1) = \Nm(1-\pi)^{m_{\pi}}.$$\end{Prop}

\begin{section}{\fontsize{12}{12}\selectfont Type $\IV(1,d)$}\end{section}

For an absolutely simple abelian variety $A$ over any field, the endomorphism algebra $\End^0(A)$ has one of the following structures according to Albert's classification:

\noindent \textbf{Type $\mr{I}$:} A totally real field.\\
\noindent \textbf{Type $\mr{II}$:} A totally definite quaternion algebra central over a totally real field.\\
\noindent \textbf{Type $\mr{III}$}: A totally indefinite quaternion algebra central over a totally real field.\\
\noindent \textbf{Type $\IV$:} A division algebra central over a CM field equipped with an involution of the second kind.

\begin{Def} \normalfont We say an absolutely simple abelian $A$ variety is of type $\IV(e,d)$ if it fulfills the following conditions:

\noindent - $D:=\End^0(A)$ is a central division algebra of  dimension $d^2$ over a CM field $K$ of degree $2e$ equipped with an involution of the second kind.\\
- The primes of $K$ ramified in $D$ lie over the same rational prime.\end{Def}

The symbol $\IV$ in the notation represents the nature of the endomorphism algebra according to Albert's classification. We impose the second condition in the definition because by Honda-Tate theory, the endomorphism algebra of a simple abelian variety over a finite field of characteristic $p$ is unramified away from the set of primes lying over $p$. The (Rosati) involution of the second kind is a result of the polarization carried by the abelian variety. Note that for any prime $\mfp$ of $K$ and its complex conjugate $\ov{\mfp}$, we have $$\mr{inv}_{\mfp}(D)+\mr{inv}_{\ov{\mfp}}(D) = 0\in \bq/\bz.$$ By Honda-Tate theory, the dimension of the simple abelian variety is $\frac{1}{2}[K:\bq][D:K]^{1/2} = ed$. In this section, we study abelian varieties of type $\IV(1,d)$ over finite fields and explore their possible use in isogeny-based cryptography. We will need the following lemma.

\begin{Lem} Let $B$ be an absolutely simple abelian variety over a finite field of characteristic $p$. If there exists a prime $\mfp$ of $K$ lying over $p$ such that the decomposition group $D_{\mfp}$ is normal in $\GalWB$, the prime $p$ splits completely in $\GalWB$.\end{Lem}

\begin{prf} This is lemma 5.3 of [Th17]\end{prf}

\begin{Prop} Let $B_1$, $B_2$ be abelian varieties of type $\IV(1,d)$ over a finite field $\bFq$. Then the following are equivalent:

\noindent - $B_1$ is isogenous to $B_2$.\\
- $\End^0(B_1)\cong \End^0(B_2)$.\end{Prop}

\begin{prf} Since isogenous abelian varieties have the same endomorphism algebra, it suffices to prove the converse. Let $\pi_{_1}$, $\pi_{_2}$ be the associated Weil numbers. Now, if $\End^0(B_1)\cong \End^0(B_2)$, then in particular, the centers $\bq(\pi_{_1})$ and $\bq(\pi_{_2})$ coincide. Let $p$ be the characteristic of $\bFq$. Since $\bq(\pi_{_1})/\bq$ is abelian, $p$ splits (completely) in $\bq(\pi_{_1})/\bq$. Let $\mfp$, $\ov{\mfp}$ be the primes of $\bq(\pi_{_1})$ lying over $p$ and let $\pi_{_1}\mc{O}_{\bq(\pi_{_1})} = \mfp^i\ov{\mfp}^j$.

Since the Hasse invariants coincide, it follows that $\pi_{_1}\pi_{_2}^{-1}$ is a unit at all non-archimedean places and hence, is a root of unity. Hence, $B_1$ and $B_2$ are twists of each other.\end{prf}

\noindent \textbf{Remark.} Note that this is the only case other than that of elliptic curves where the endomorphism algebra determines the isogeny class. For any other division algebra $D$ of type $\IV(e,d)$ with $e\geq 2$ it is fairly easy to construct a pair of non-isogenous abelian varieties with endomorphism algebra $D$.

Let $B$ be an abelian variety of type $\IV(1,d)$ over a finite field $\bFq$. Let $\mc{O}_D$ be a fixed maximal order in $D:=\End^0(B)$ and let $\Iso(B)$ denote the set of abelian varieties over $\bFq$ isogenous to $B$. The map \vspace{-0.15cm}$$\Cl(\mc{O}_D)\times \Iso(B)\lra \Iso(B),\;\;\; (I, B')\mapsto IB'$$ makes $\Iso(B)$ into a $\Cl(D)$-torsor. Hence, $\Iso(B)$ has cardinality $\# \Cl(\mc{O}_D)$. Unlike the ordinary case, the class group $\Cl(\mc{O}_D)$ is not a commutative group. In fact, $\Cl(\mc{O}_D)$ has the structure of a pointed set rather than a group. This makes the resulting cryptosystem more resistant to the quantum attacks similar to those outlined in [CJS]. 

\begin{subsection}{\fontsize{11}{11}\selectfont Newton Polygons}\end{subsection}


Let $B$ be an abelian variety over an algebraically closed field $k$ of characteristic $p>0$. The group scheme $B[p^{\infty}]$ is a $p$-divisible group of rank $\leq \dim B$. Let $D(B[p^{\infty}])$ denote the Dieudonne module. Then $D(B[p^{\infty}])\otimes_k W(k)[\frac{1}{p}]$ is a direct sum of pure isocrystals by the Dieudonne-Manin classification theorem. Let $\lamb_1<\cdots<\lamb_r$ be the distinct slopes and let $m_i$ denote the multiplicity of $\lamb_i$. The sequence $m_1\times \lamb_1,\cdots, m_r\times \lamb_r$ is called the \textit{Newton polygon} of $B$.\\

\begin{Def} \normalfont A Newton polygon is \textit{admissible} if it fulfills the following conditions:

\noindent 1. The breakpoints are integral, meaning that for any slope $\lambda$ of multiplicity $m_{\lambda}$, we have $m_{\lambda}\lambda \in \bz$.\\ 
\noindent 2. The polygon is \textit{symmetric}, meaning that each slope $\lambda$, the slopes $\lambda$ and $1-\lambda$ have the same multiplicity.\end{Def}

Let $\pi$ be a Weil $q$-integer and let $B_{\pi}$ be the corresponding simple abelian variety over $\bFq$. Then the Newton slopes of $B_{\pi}$ are given by $\{v(\pi)/v(q)\}_v$ where $v$ runs through the places of $\bq(\pi)$ lying over $p$. In particular, the Newton polygon is symmetric and hence, all slopes lie in the interval $[0,1]$. A far more subtle fact is that the converse is also true. This was formerly known as Manin's conjecture until proven by Serre. We refer the reader to [Tat69] for the proof.

\begin{Thm} $\mr{(Serre)}$ The Newton polygon of an abelian variety over a finite field is admissible. Conversely, any admissible polygon occurs as the Newton polygon of some abelian variety over a finite field of any prescribed characteristic.\end{Thm}

\begin{Prop} Let $p$ be a rational prime and $K$ an imaginary quadratic field in which $p$ splits. Let $D$ be a central division algebra over $K$ that is unramified away from the places lying over $p$. Then there exists an abelian variety $B$ over a finite extension $\bFq$ of $\bFp$ such that $\End^0(B)\cong D$. Furthermore, there exists an abelian variety $\wti{B}$ over the field $\bFq$ with the endomorphism ring $\End(B)$ a maximal order in $\End^0(B)$.\end{Prop}

\begin{prf} Let $\mfp$, $\ov{\mfp}$ be the primes of $K$ lying over $\bq$. Since $D$ is split at all primes other than $\mfp,\ov{\mfp}$, we have \vspace{-0.15cm}$$\inv_\mfp(D)+\inv_{\ov{\mfp}}(D) = \sum_v \inv_v(D) = 0 \in \bq/\bz.$$ Furthermore, the least common denominator of $\{\inv_\mfp(D), \inv_{\ov{\mfp}}(D) \}$ is $d$ and hence, $\inv_\mfp(D) = \frac{j}{d}$ for some $j<d$ prime to $d$. Choose an integer $n$ such that $\mfp^n$ is principal and write $\mfp^n = (\pi)$. Now, $\pi_{_1}:=\pi^j\ov{\pi}^{d-j}$ is a Weil $p^{dn}$-integer and the abelian variety $B_{\pi_{_1}}$ corresponding to $\pi_{_1}$ has Newton slopes $$\frac{v_{\mfp}(\pi_{_1})}{v_{\mfp}(p^{dn})} = \frac{j}{d},\;\;\frac{v_{\mfp}(\ov{\pi}_{_1})}{v_{\mfp}(p^{dn})} = \frac{d-j}{d}.$$ Hence, the endomorphism algebra of $B_{\pi_{_1}}$ is isomorphic to $D$.

Let $\Lamb$ be a maximal order of $D$ that contains the order $\End(B)$. By ([Yu11], Theorem 1.3), there exists an abelian variety $\wti{B}$ over the same field of definition as $B$ such that $\End(\wti{B}) = \Lamb$ and there is a universal isogeny $\phi:\wti{B}\lra B$ such that any isogeny $B'\lra B$ factors through $\phi$.\end{prf}

\noindent \textbf{Remark} Note that an abelian variety of type $\IV(1,1)$ is just an ordinary elliptic curve. There does not exist an abelian variety of type $\IV(1,2)$. Any such abelian variety would have Newton polygon $4\times 1/2$ which only occurs for a second power of a supersingular elliptic curve. 

\begin{Prop} Let $B_{\pi}$ be a simple abelian variety over $\bFP$ of dimension $d$ corresponding to a Weil number $\pi$. The following are equivalent:

\noindent $(1)$ $B_{\pi}$ is of type $\IV(1,d)$.\\
$(2)$ $\bq(\pi)/\bq$ is abelian and $B_{\pi}$ has Newton polygon \vspace{-0.1cm}$$d\times \frac{j}{d},\;\; d\times \frac{d-j}{d}$$ for some integer $j$ prime to $d$.\end{Prop}

\begin{prf} $(1)\Rightarrow (2)$: $\bq(\pi)/\bq$ is imaginary quadratic and in particular, is an abelian extension. So $p$ splits in $\bq(\pi)$. The dimension of $\End^0(B_{\pi})$ is the least common multiple of its Hasse invariants. Furthermore, since the primes of $\bq(\pi)$ lying over $p$ have local degree one over $\bq$, the set of Newton slopes of $B_{\pi}$ coincides with the set of the Hasse invariants of $\End^0(B_{\pi})$. Hence, the least common denominator of the Newton slopes is $d$, which implies $(2)$.

\noindent $(2)\Rightarrow (1)$: Since $\bq(\pi)/\bq$ is abelian, $p$ splits completely in $\bq(\pi)$. Hence, the dimension of $\End^0(B_{\pi})$ over $\bq(\pi)$ is $d^2$. Now, \vspace{-0.1cm}$$2d = [\bq(\pi):\bq][\End^0(B_{\pi}):\bq(\pi)]^{1/2} = d[\bq(\pi):\bq]$$ and hence, $[\bq(\pi):\bq] = 2$.\end{prf}

\begin{Corr} Any abelian variety $B$ of type $\IV(1,d)$ over a finite field of characteristic $p$ has $p$-rank zero.\end{Corr}

\begin{prf} The $p$-rank of an abelian variety over a field of characteristic $p$ is the multiplicity of the slopes $0$ and $1$. Since these slopes do not occur in the Newton polygon of $B$, it follows that $B[p] = \{0_B\}$.\end{prf}

Thus, the Frobenius and the Verschiebung are purely inseparable morphisms. It is well-known that every abelian variety is isogenous to a principally polarized abelian variety over the algebraic closure. But in this case, we do not need to pass to a larger extension to have a ppav in the isogeny class.

\begin{Prop} Any abelian variety $B$ of type $\IV(1,d)$ over a finite field $k$ is isogenous to a principally polarized abelian variety over the same field $k$.\end{Prop}
\begin{prf} Let $p$ be the characteristic $k$ and let $\phi:B\lra \what{B}$ be a polarization from $B$ to its dual. Since $B$ has $p$-rank zero, we may assume the polarization is separable. By [Mil], Proposition ???, the isogeny class of $B$ over $k$ contains a principally polarized abelian variety.\end{prf}

\begin{Prop} Let the notations be as in the last proposition and let $\mfl$ be a prime of the imaginary quadratic field $K:=\Cent(\End^0(B))$ that does not divide $\charac(k)$. Then the degree $[k(B[\mfl]):k]$ is the smallest integer $N$ such that $\pi_B^{N}\equiv 1\pmod{\mfl}$ in $\mc{O}_K$.\end{Prop}

\begin{subsection}{\fontsize{11}{11}\selectfont The minimum field of definition}\end{subsection}

\begin{Def} \normalfont For an abelian variety $A$ over a field $F$, we say $A$ has a \textit{model} over a subfield $F_0$ if there exists an abelian  variety $A_0$ over $F_0$ such that the base change $A_0\times_{F_0} F =_{\isog} A$.\end{Def}

\noindent We show that in our setting, the class number of the center gives an upper bound on the minimum field of definition.

\begin{Prop} Let $B$ be an abelian variety of type $\IV(1,d)$ with endomorphism algebra $D$ and let $K$ be the center of $D$. Let $h$ be the class number of $K$. Then $B$ has a model over $\bFq$ where $q := p^h$.\end{Prop}

\begin{prf} Let $\pi$ be the Weil number corresponding to $B$ and let $\mfp$, $\ov{\mfp}$ be the primes of $K$ lying over $p$. Write $\pi_{_0}\mc{O}_K= \mfp^h$. Then $\pi_{_0}$ is a Weil $p^h$-integer and hence, corresponds to an abelian variety $B_{\pi_{_0}}$ of dimension $d$ over $\bFq$. Furthermore, $B_{\pi_{_0}}$ has endomorphism algebra $D$ and hence, is isomorphic to $B$ over some finite extension.\end{prf}

Let $B$ be an abelian variety with endomorphism algebra $D$. Write $\Lamb = \End(B)$ and let $I$ be a left ideal of $\Lamb$. We associate to $I$ an isogeny $\phi_I$ as follows. Consider the finite group scheme \vspace{-0.1cm}$$C_I := \bigcap\limits_{i\in I} \ker(i: B\lra B) .$$ This is a finite subgroup scheme of $B$ and hence, yields an isogeny \vspace{-0.1cm}$$\phi_I:B\lra B^{(I)}:=  B/C_I.$$ We call this the \textit{ideal isogeny} associated to $I$. In the case of supersingular elliptic curves, it is well-known that every isogeny is an ideal isogeny arising from some left ideal of the endomorphism ring. This is a consequence of the fact that the center of the endomorphism algebra of a supersingular elliptic curve is $\bq$. On the other hand, for any simple non-supersingular abelian variety over a finite field, the center of the endomorphism algebra is a CM field. Consequently, any isogeny of a degree that is not a norm of the CM field does not arise as an ideal isogeny.

Conversely, for any isogeny $\phi:B\lra B'$, the annihilator \vspace{0.1cm}$$\mr{ann}(\ker(\phi)) := \{\al\in \End(B):\; \al(\ker(\phi)) = \{0_B \} \}$$ is a left ideal of $\End(B)$.


\begin{Prop} Let $D$ be a division algebra of type $\IV(1,d)$ for some odd integer $d\geq 3$. There exists an abelian variety $A$ with complex multiplication over a number field $F$ and a prime $v$ of good reduction such that $A_v$ is of type $\IV(1,d)$.\end{Prop}

\begin{prf} Let $K$ be the center of $D$ and let $p$ be the characteristic of the primes of $K$ ramified in $D$. Choose a degree $d$ cyclic extension $L_0$ over $\bq$ such that the $p$ is inert in $L_0$. Then $L:=KL_0$ is cyclic of degree $d$ over $K$ and the primes of $K$ lying over $p$ are inert in $L$. Let $X$ be the abelian variety $\bc^d/\mc{O}_L$ over $\bc$. By a theorem of Shimura, $X$ has a model $A$ over some number field $F$. Since $A$ has complex multiplication, it has potential good reduction everywhere. So, replacing $F$ by a finite extension if necessary, we may assume $A$ has good reduction everywhere. Let $v$ be a prime of $F$ lying over $p$. Then each component of $A_v[p^{\infty}]$ is isoclinic of slope $j/d$ for some integer $j$ prime to $d$. Thus, $\End^0(A_v)$ is a central division algebra of dimension $d^2$ over $K$.\end{prf}

\begin{Def} \normalfont An isogeny $\phi:B_1\lra B_2$ of abelian varieties over a finite field is said to be \textit{cyclic} if the degree is prime to the characteristic of the field and the kernel $\ker(\phi)$ is a cyclic subgroup scheme of $\ker(\deg \phi)$.\end{Def}

Note that when the abelian varieties $B_i$ are elliptic curves, this is equivalent to the isogeny $\phi$ \text{not} being of the form $[n]\circ \phi'$ for some isogeny $\phi'$. Thus, it coincides with the notion of a \textit{primitive isogeny} studied in [Koh96].

\begin{Example}\normalfont Let $E$ be the elliptic curve $y^2:= x^3 + 2$ over the finite field $\bF _{13}$. For any integer $n$ prime to $7$, the endomorphism \vspace{-0.1cm}$$[n]:E\lra E,\;P\mapsto nP$$ has kernel isomorphic to $(\bz/n\bz)^2$ and hence, is not cyclic. On the other hand, the endomorphism corresponding to the algebraic integer $2+\sqrt{-1}$ has kernel isomorphic to $\bz/5\bz$ and hence, is a cyclic isogeny.\end{Example}

In the case of supersingular elliptic curves, it is well-known that every isogeny arises from a left ideal of the endomorphism ring. Although this behavior is not shared by abelian varieties of type $\IV(1,d)$, they exhibit some similar properties. To this end, we will need the following notion.

\begin{Def} \normalfont Let $\phi:A_1\lra A_2$ be an isogeny of abelian varieties of type $\IV(1,d)$ and let $K$ be the center of the endomorphism algebra $\End^0(A_1)$. We say the isogeny $\phi$ is a \textit{normed isogeny} if the set \vspace{-0.1cm}$$I(\ker(\phi)):=\{\al\in \mc{O}_K: \al(\ker(\phi)) = 0 \} $$ is stable under the involution on $\End^0(A_1)$.\end{Def}

\noindent (This non-standard terminology is a placeholder until we come across a better term in the existing literature).

\begin{Prop} The following are equivalent:

\noindent $(1)$ $l$ is a norm of an integral element of $K$.\\
\noindent $(2)$ $l$ splits completely in $K$ and the primes of $K$ lying over $l$ are principal.\\
\noindent $(3)$ $l$ splits completely in the Hilbert class field $\mr{H}(L)$ of $L$.\end{Prop}

\begin{prf} Note that since $K/\bq$ is assumed to be a Galois extension, the primes of $K$ lying over $l$ are Galois conjugates and hence, the corresponding ideal classes have the same orders as elements of the ideal clas group.

\noindent $(1)\Leftrightarrow (2)$: Suppose $l = \Nm_{K/\bq}(\al)$ for some $\al\in \mc{O}_K$. Then $l\bz = \Nm_{K/\bq}(\al\mc{O}_K)$ and hence, $\al\mc{O}_K$ is a prime of $K$ lying over $l$. Conversely, if a prime $\mfl$ of $K$ lying over $l$ is principal with $\mfl = \be\mc{O}_K$, then $\Nm_{K/\bq}(\be) = l$.

\noindent$(2)\Leftrightarrow (3)$: Let $\mfl$ be a prime of $K$ lying over $\bq$. By the principal ideal theorem, $\mfl$ is a principal ideal if and only if it splits completely in the Hilbert class field $\mr{H}(L)$. If $l$ splits completely in $K/\bq$, this is equivalent to $l$ splitting completely in $\mr{H}(L)/\bq$.\end{prf}

\begin{Prop} Let $B$ be an abelian variety of type $\IV(1,d)$ with the endomorphism ring $\End(B)$ a maximal order in the endomorphism algebra $\End^0(B)$. Let $I$ be a left ideal of $\Lamb$ and let $\phi_I:B\lra B/C_I$ be the ideal isogeny associated to $I$. Then $\phi_I$ is a normed isogeny. Furthermore, the endomorphism ring of $B' := B/C_I$ is a maximal order in $D$.\end{Prop}

\begin{prf} Since $O_l(I) = \Lamb$ is assumed to be a maximal order, so is the right order $O_r(I)$ of $I$. The degree of $\phi_I$ is given by $\Nm_{K/\bq}(\Nr(I))$ and hence, $\phi_I$ is a normed isogeny. Choose an integer $n$ such that $nO_r(I)\sub \Lamb$ and let $\be\in O_r(I)$. Now $nI\be\sub nI$ and hence, $\ker(nI)\sub \ker(\phi_I\circ n\be)$. So $\phi_I\circ n\be$ factors through $n\phi_I$ and hence, $\be$ induces an endomorphism of the abelian variety $B^{(I)}:=\phi_I(B)$.	Thus, we have an injective homomorphism $O_r(I)\hra \End(B^{(I)})$. Since $O_r(I)$ is a maximal order in $D$, this inclusion is an equality, which completes the proof.\end{prf}

Later, we will prove the converse.

\begin{Prop} The following are equivalent:

\noindent $(1)$ $\psi_2$ lies in the ring of integers of the field $K(\psi_1)$.\\
$(2)$ $\psi_1$ commutes with $\psi_2$.\end{Prop}

\begin{prf} The field $K(\psi)$ is of degree $d$ over $K$ and hence, is a maximal subfield of $D$. By the double centralizer theorem, the field $K(\psi)$ is its own centralizer in $D$. Hence, $\psi_2\in K(\psi)$. Since $\psi_2$ is an isogeny, it is integral and hence, lies in the ring of integers of $K(\psi_1)$.\end{prf}

\noindent The following decomposition explains our motivation for studying Galois endomorphisms.  

\begin{Prop} Let $\psi_1,\cdots,\psi_r$ be pairwise commutative endomorphisms of $B$ of degree prime to the characteristic of the field of definition. Then there exists a maximal subfield $L$ of $D$ such that $\psi_1,\cdots,\psi_r$ have embeddings in $\mc{O}_L$. Furthermore, if for each pair $i,j $, the elements $\psi_i$ and $\psi_j$ generate relatively prime ideals of $K(\psi_1,\psi_2)$, then \vspace{-0.1cm}$$B[\prod\limits_{i=1}^r \psi_i]= \bigoplus\limits_{i=1}^r B[\psi_i].$$\end{Prop}

\begin{prf} We show this for the case where $r=2$. The general case then follows by induction.

Since $K(\psi_1,\psi_2)$ is a subfield of $D$, it is contained in some maximal subfield $L$ of degree $d$ over $K$. Now, $\psi_1$ and $\psi_2$ generate relatively prime ideals in $R:=\End(B)\cap \bz[\psi_1,\psi_2]$. Hence, there exist $\al_1,\al_2\in R$ such that $\al_1\psi_1 + \al_2\psi_2 = 1$ in $R$. Now, for any point $P$ of $B$, if $\psi_1(P) = \psi_2(P) = 0$, then \vspace{-0.1cm}$$P = \al_1\psi_1(P) + \al_2\psi_2(P) = 0_B.$$ Thus, $B[\psi_1]\cap B[\psi_2] = \{ 0_B \}$. 

It remains to show that $B[\psi_1\psi_2]\subseteq B[\psi_1]\oplus B[\psi_2]$. Let $P\in B[\psi_1\psi_2]$. Then $\psi_1(P)\in B[\psi_2]$ and $\psi_2(P)\in B[\psi_1]$. In particular, $\al_1\psi_1(P)\in B[\psi_2]$, $\al_2\psi_2(P)\in B[\psi_1]$ and hence, \vspace{-0.1cm}$$P = \al_2\psi_2(P)+\al_1\psi_1(P)\in B[\psi_1]\oplus B[\psi_2],$$ which completes the proof.\end{prf}

\noindent \textbf{Remark.} Note that this does not necessarily hold when the endomorphisms $\psi_1$ and $\psi_2$ do not commute. We will need the following notion.



\begin{Def} Let $A$ be an absolutely simple abelian variety and let $F$ be the center of its endomorphism algebra $\End^0(A)$. We say an endomorphism $\phi:A\lra A$ is \textit{Galois} if there exists a maximal subfield $L$ of the division algebra $\End^0(A)$ such that $L/F$ is a Galois extension and $F(\phi)$ is a subfield of $L$.\end{Def}

\noindent This definition extends naturally to quasi-endomorphisms of the abelian variety. Note that in the case where $\End^0(A)$ is a quaternion algebra, this notion is clearly vacuous since every quadratic extension is Galois. On the other hand, a division algebra of dimension larger than $4$ over a number field will have infinitely many maximal subfields that are non-Galois over $F$.

\noindent (This nonstandard terminology is a placeholder until we come across a better term in the existing literature).

\begin{Prop} If $\phi_1,\phi_2$ are Galois endomorphisms of a simple abelian variety $A$ that commute with each other, then the endomorphisms $\phi_1+\phi_2$ and $\phi_1\phi_2$ are Galois as well.\end{Prop}
\begin{prf} We omit the proof since it is straightforward.\end{prf}

\begin{Prop} Let $\psi$ be a Galois endomorphism of $B$ with reduced norm $\Psi$. Suppose the ideals in $\mc{O}_{F(\psi)}$ generated by $\psi$ and $\Psi\psi^{- d/[F(\psi):F]}$ are relatively prime. Then we have a decomposition \vspace{-0.15cm}$$B[\Psi] = \bigoplus\limits_{j=1}^r B[\psi_j] $$ where the $\psi_j$ are the Galois conjugates of $\psi$ over $F$.\end{Prop}

\begin{prf} By assumption, there exists a maximal subfield $L$ of $D$ such that $L/F$ is a Galois extension and $F(\psi)\sub L$. Hence, the Galois closure of $F(\psi)$ over $F$ has an embedding in $L$. Let $\psi_1,\cdots,\psi_r$ be the distinct conjugates of $\psi$ over $F$, where $r = [F(\psi):F]$. Then \vspace{-0.1cm}$$\Psi := \Nr(\psi) = \Nm_{L/\bq}(\psi) = \prod\limits_{i=1}^{r}\psi_i^{[L:F(\psi)]}.$$ Now, by assumption, the ideals of $\mc{O}_L$ generated by the elements $\psi_i^{[L:F(\psi)]}$ are pairwise co-prime. Since the $\psi_i$ are elements of the same maximal subfield $L$, they are pairwise commutative. Hence, by the preceding proposition, we have the decomposition \vspace{-0.1cm}$$B[\Psi] = \bigoplus\limits_{j=1}^r B[\psi_j^{[L:F(\psi)]}],$$ which completes the proof. \end{prf}

\begin{Prop} Let $\al:B\lra B$ be a Galois endomorphism of reduced norm $a\in \mc{O}_K$. Then $\al$ is cyclic if and only if it satisfies all of the following conditions:

\noindent $1$. $L:=K(\al)$ is a Galois extension degree of $d$ over $K$.\\
$2$. For any $\sigma\in \Gal(L/K)$ other than the identity, the ideals $\al\mc{O}_L$ and $\sigma(\al)\mc{O}_L$ are relatively prime.\\
$3$. The primes of $L$ dividing $\al\mc{O}_L$ are of local degree one over $K$ and for any prime $\mfp$ of $K$ that divides $\Nr(\al)\mc{O}_K$, precisely one prime of $L$ lying over $\mfp$ divides $\al\mc{O}_L$.\end{Prop}

\begin{prf} Set $a_0:= \Nm_{L/K}(\al)$ and $a := a_0a_0^*$. Then $B[a_0]\cong B[a_0^*]\cong (\bz/a\bz)^d$. Furthermore, \vspace{-0.1cm}$$B[a] = \sum\limits_{\sigma\in \Gal(L/K)} B[\sigma(\al)]$$ and the group structure of $B[\sigma(\al)]$ is independent of $\sigma$. Hence, each torsion subgroup scheme $B[\sigma(\al)]$ is cyclic if and only if $$B[\sigma(\al)]\cap B[a\sigma(\al)^{-1}] = \{0_B\} \;\forall\; \sigma\in \Gal(L/K).$$ But this is equivalent to the ideals $\al\mc{O}_L$ and $a\sigma(\al)^{-1}\mc{O}_L$ being relatively prime.\end{prf}




\begin{subsection}{\fontsize{11}{11}\selectfont The set of left ideal classes}\end{subsection}

Let $K$ be an imaginary quadratic field and let $p$ be a prime that splits in $K$. Let $D$ be a division algebra of dimension $d^2$ ramified only at the primes lying over $p$. So $D$ has Hasse invariants $\{j/d, -j/d\}$ at the places lying over $p$ for some integer $j$ prime to $d$. We fix a maximal order $\Lamb$ in the division algebra $D$.

\begin{Def} \normalfont Two left ideals $I$, $J$ of $\Lamb$ are equivalent if $J = Ia$ for some $a\in D^{\times}$. The number $\mr{h}(\Lamb)$ of equivalence classes of left $\Lamb$-ideals is called the class number of $\Lamb$.\end{Def}

Unlike the commutative case, the left ideals do not form a group under multiplication. It is simply a pointed set with the distinguished element corresponding to the class of principal left ideals. This lack of a group structure makes it less susceptible to quantum attacks similar to those outlined in [CJS].

If $\Lamb$ is a maximal order in a division algebra of type $\IV(1,d)$, the class number is finite, by the next proposition. Furthermore, it is well-known that any two maximal orders are conjugate by some element of the idele group $$D_{\mb{A}}^{\times} = \{(a_v)_v \in \prod\limits_v D_v^{\times}\;|\; a_v \in \Lamb_v \text{ for all but finitely many } v \}.$$ Hence, the class number is independent of the choice of the maximal order. Furthermore, we have a natural norm map defined as follows:\vspace{-0.15cm}$$D_{\mb{A}}^{\times}\lra \bc^{\times},\;\;\; a\mapsto \prod\limits_{v}\big| \Nm(a_v) \big|_v.$$ We denote the kernel of this norm map by $D_{\mb{A}}^{(1)}$. By the product formula, the group $D_K^{\times}$ embeds diagonally into $D_{\mb{A}}^{(1)}$.

\begin{Prop} Let $\Lamb$ be a maximal order in $D$. The class number of $\Lamb$ is finite.\end{Prop}
\begin{prf} If $v$ is any place of $K$ prime to $p$, then $D$ splits at $v$. Choosing the isomorphism $D_v\cong \Mat_{d}(K_v)$ appropriately, we may assume that $\Lamb_v = \Mat_{d}(\mc{O}_{K_v})$. So every left ideal $I$ of $\Lamb$ is principal at the places where $D$ splits (which in this case is all places prime to $p$). For all but finitely many places, we have $I_v = \Lamb_v$, meaning $a_v\in \Lamb_v^{\times}$. Hence, there exists an element $b\in D_{\mb{A}}^{(1)}$ such that $b_v = a_v$ at all but finitely many places. Hence, $D_{\mb{A}}^{(1)}$ acts transitively (from the left) on the set of left ideal classes. 

The actions of $D_K^{\times}$ and $\Lamb_{\mb{A}}^{\times}$ on the set of ideal classes are trivial and hence, the class number equals the number of double cosets in $\Lamb_{\mb{A}}^{\times}\setminus D_{\mb{A}}^{(1)}/D_K^{\times}$. Since $D^{\times}$ is a connected reductive group, it follows that that $D_{\mb{A}}^{(1)}/D_K^{\times}$ is compact and $\Lamb_{\mb{A}}^{\times}$ is open in $D_{\mb{A}}^{(1)}$. Hence, the class number is finite.\end{prf}

\begin{Def} \normalfont Let $D$ be a division algebra central over a global field $F$. Let $M$ be a full $\mc{O}_F$-lattice in $D$. The \textit{left} and \textit{right orders} of $M$ are defined as follows:\vspace{-0.15cm} $$O_l(M):=\{a\in \mc{O}_K:\; Ma\sub M \},\;\;\;O_r(M):=\{a\in D:\; aM\sub M \}.$$\end{Def}

\noindent The following theorem is fundamental to the theory of orders in central simple algebras.

\begin{Thm} For any full $\mc{O}_F$-lattice $M$, the order $O_l(M)$ is a maximal order if and only if $O_r(M)$ is.\end{Thm}
\begin{prf} See ([Rei75], Chapter 3).\end{prf}

\noindent An immediate implication is the following.

\begin{Corr} Let $\Lamb$ be a maximal left order in $D$ and let $I$ be a left ideal of $\Lamb$. Then $O_r(I)$ is a maximal order in $D$.\end{Corr}
\begin{prf} Since $I$ is a left ideal, clearly $\Lamb\sub O_l(I)$ and since $\Lamb$ is a maximal order, we have $\Lamb = O_l(I)$. By the preceding theorem, $O_r(I)$ is a maximal order.\end{prf}

\begin{Prop} Let $\Lamb$ be a maximal order in $D$ and let $\{I_1,\cdots,I_n \}$ be a set of left $\Lamb$-ideals that represent the distinct left ideal classes of $\Lamb$. Then each conjugacy class of maximal orders in $D$ is represented in the set of right orders $\{O_r(I_1),\cdots, O_r(I_n) \}.$\end{Prop}

\begin{prf} Let $\Lamb'$ be any maximal ideal of $D$. For any place $v$ of $K$, the orders $\Lamb_v$ and $\Lamb_v '$ are maximal orders in $\mc{O}_{K_v}$. Hence, there exists $\al_v\in D_v^{\times}$ such that $\al_v^{-1}\Lamb_v\al_v = \Lamb_v$. Now, there exists $\wti{\al}\in D_{\mb{A}}^{(1)}$ such that $\wti{\al}_v = \al_v$. Set $I:= \Lamb\wti{\al}$. Then $I$ is a left ideal of $\Lamb$ with right order $\Lamb'$. Since the ideals $\{I_1,\cdots,I_n \}$ represent all of the ideal classes of $\Lamb$, there exists $a\in D^{\times}$ and an index $i$ such that $J=I_ia$. Thus, $O_r(I_i) = a\Lamb'a^{-1}$.\end{prf}

\begin{Prop} Let $D$ be a division algebra of type $\IV(1,d)$ with center $K$, $\Lamb$ a maximal order in $D$ and $\mfl$ a prime of $K$ with characteristic $l$. There exists an element $\al\in \Lamb$ such that any maximal left ideal of $\Lamb$ containing $\mfl$ is of the form $J= (\Lamb l+\Lamb \al)u$ for some $u\in D^{\times}$ with reduced norm $\Nr(u) = \pm 1$.\end{Prop}

\begin{prf} We may choose a Galois extension $L/K$ of degree $d$ such that:

\noindent - $L$ is inert at the primes of $K$ that $D$ is ramified at.\\
\noindent - $\mfl$ splits completely in $L$.\\
\noindent - $L$ is linearly disjoint from the spinor class field of $K$.

The first condition ensures that $L$ has an embedding in $D$ and the last condition ensures that the ring of integers $\mc{O}_L$ has an embedding in every maximal order of $D$. We fix an embedding $\mc{O}_L\hra \Lamb$. Now, let $\wti{\mfl}$ be a prime of $L$ lying over $\mfl$ and choose an element $\al$ of $\wti{\mfl}$ such that $\wti{\mfl} =  l\mc{O}_L+\al\mc{O}_L$ as an ideal in $\mc{O}_L$. Let $\mf{m}$ be a maximal left ideal of $\Lamb$ containing the left $\Lamb$-ideal $\Lamb\wti{\mfl}:=\{\sum\limits_{x\in \Lamb, \be\in \wti{\mfl}} x\be \}$. Then $\mfl:= \mf{m}\cap \mc{O}_K$ is a prime in $K$ and $\Nr(\mf{m}) = \mfl$ (c.f.  [Rei75]). On the other hand, $\Nm_{L/K}(\wti{\mfl}) = \mfl$ and hence, the inclusion $\Nr(\Lamb\wti{\mfl})\sub \mfl = \Nr(\mfl)$ is an equality. Thus, $\Lamb\wti{\mfl}$ is a maximal left ideal that contains $\mfl$ and hence, the inclusion $\Lamb\wti{\mfl}\sub \mf{m}$ is an equality.

Now let $J$ be any maximal left ideal of $\Lamb$ containing $\mfl$. Then $\Nr(J) = \mfl =\Nr(\mf{m}) $ (c.f.  [Rei75]) and since $D$ fulfills Eichler's condition, Eichler's refinement of the Hass-Mass--Schilling theorem implies that $J = \mf{m}u$ for some $u\in D^{\times}$. Thus, $J = \Lamb l u + \Lamb \al u$. Furthermore, $\Nr(u)$ is a unit in $\mc{O}_K$ and since $K$ is an imaginary quadratic field, $\Nr(u) = \pm 1$.\end{prf}

\begin{Prop} Let $B$ be an abelian variety of type $\IV(1,d)$ with $\End(B)$ a maximal order in $D:= \End^0(B)$. Let $J$ be a maximal left ideal of $\End(B)$ with the prime $\mfl:= J\cap K$ a prime of local degree one over $\bq$. Then the isogeny $\phi_J$ is cyclic of prime degree.\end{Prop}

\begin{prf} We write $\Lamb:=\End(B)$ for brevity. We choose a degree $d$ extension $L/K$ as in the preceding proposition and fix an embedding $\mc{O}_L\hra \Lamb$. Now, let $\wti{\mfl}$ be a prime of $L$ lying over $\mfl$ and choose an element $\al_1$ of $\wti{\mfl}$ such that $\wti{\mfl} = l\mc{O}_L+\al\mc{O}_L$ as an ideal in $\mc{O}_L$.  Let $\al_1,\cdots,\al_d$ be the distinct Galois conjugates of $\al_1$ in the fixed embedding of $\mc{O}_L$ in $\Lamb$ and write $\wti{\mfl}_i:=l\mc{O}_L+\al_i\mc{O}_L$, $i=1,\cdots,d$. By the argument in the last proposition, the left ideals $\mf{m}_i:=\Lamb \wti{\mfl}_i$ are maximal left ideals in $\Lamb$. 

Note that the elements $\al_1,\cdots,\al_d$ are pairwise commutative by construction and hence, for each index $i$, the product \vspace{-0.1cm}$$\prod\limits_{i\neq j}^d \al_i \in \bigcap\limits_{i\neq j}^d \mf{m}_i,\;\;\notin \mf{m}_j .$$ Thus, the left $\Lamb$-ideals $\mf{m}_j$ and $\bigcap\limits_{i\neq j} \mf{m}_i$ are co-maximal and hence, \vspace{-0.1cm}$$\ker(\mfl) = \bigoplus\limits_{i=1}^d \ker(\mf{m}_i),$$ which implies that $\ker(\phi_{\mf{m}_i}) \cong \bz/l\bz$ where $l$ is the rational prime lying under $\mfl$.

Now let $\mf{m}'$ be any maximal left ideal of $\Lamb$ with $\mf{m}'\cap K = \mfl$. Then Eichlers's refinement of the Hasse-Mass-Schilling theorem implies that $\mf{m}' = \mf{m}_1u$ for some $u\in D^{\times}$ with $\Nr(u)=\pm 1$. Hence, $\ker(\phi_{\mf{m}'}) \cong \ker(\phi_{\mf{m}_1})$.\end{prf}



\begin{Prop} Let $D$ be a division algebra of type $\IV(1,d)$. For any maximal $\Lamb$ order of $D$, there exists an abelian variety $B$ over $\bFP$ with $\End(B) = \Lamb$.\end{Prop}

\begin{prf} As shown in Proposition 2.4, there exists an abelian variety $B'$ with its endomorphism ring $\Lamb'$ a maximal order in $D$. By the preceding proposition, there exists a left ideal $I'$ of $\Lamb$ such that $O_r(I') = \Lamb$. The abelian variety \vspace{-0.1cm}$$B:=B'^{(I')}=\phi_{I'}(B')$$ is isogenous to $B'$ and has endomorphism ring $\Lamb$.\end{prf}

\begin{subsection}{\fontsize{11}{11}\selectfont The image of the reduced norm}\end{subsection}

\noindent \textbf{The reduced norm.} Let $D$ be a central division algebra of dimension $d^2$ over a field $F$. Let $\al$ be an element of $D^{\times}$. Now, $\al$ commutes with $F$ and hence, $F(\al)$ is a field with the degree $[F(\al):F]$ dividing $d$. Let $L$ be a maximal subfield of $D$ such that $F(\al)\sub L$. The field norm $\Nm_{L/F}(\al)$ is called the \textit{reduced norm} of $\al$ in $D/F$. We denote the reduced norm by $\Nr_{_D}(\al)$ or simply by $\Nr(a)$ when there is no ambiguity.

Let $A$ be a simple abelian variety and let $F$ be the center of its endomorphism algebra $\End^0(A)$. For any endomorphism $\phi:A\lra A$, the degree of $\phi$ is given by \vspace{-0.1cm}$$\deg \phi = \Nm_{F/\bq}(\Nr(\al)).$$

\begin{Thm} \normalfont{(Hasse-Maas-Schilling)} Let $D$ be a division algebra central over a global field $F$. For any element $a\in F^{\times}$ that is positive at all real archimedean places that ramified in $D$, there exists an element $\al\in D^{\times}$ such that $\Nr(\al) = a$.\end{Thm}

Since we are concerned with isogenies rather than quasi-isogenies, it is necessary to consider the question of which elements of the field occur as norms as elements of $D$ integral over $\mc{O}_F$. The Hasse-Mass-Schilling theorem was refined by Eichler for the cases when the division algebra satisfies the following reasonable conditions.

\begin{Def} \normalfont A central simple algebra $A$ over a number field $F$ is said to fulfill \textit{the Eichler condition} unless both of the following properties hold:

\noindent - $F$ is totally real.\\
- $A$ is a totally definite quaternion algebra over $F$.\end{Def}  

When $A$ satisfies both conditions, we say that $A$ \textit{fails} Eichler's condition. In this article, we are concerned with the cases where $A$ is the endomorphism algebra of an abelian variety over a finite field. The following is immediate from Honda-Tate theory.

\begin{Prop} Let $B$ be an absolutely simple abelian variety over a finite field. The following are equivalent:

\noindent - $\End^0(B)$ fails Eichler's condition.\\
- $B$ is a supersingular elliptic curve.\end{Prop}

\begin{prf} If $B$ is a supersingular elliptic curve, then $\End^0(B)$ is a quaternion algebra over the totally real field $\bq$ and is ramified at the infinite place of $\bq$. In all other cases, the center of $\End^0(B)$ is a CM field and hence, $\End^0(B)$ fulfills Eichler's condition.\end{prf}

The following is a refinement of the Hasse-Maas-Schilling theorem for the case where $D$ is Eichler over $F$. 

\begin{Thm} $\mr{(Eichler)}$ Let $D$ be a division algebra central over a global field $F$ that fulfills Eichler's condition and let $\Lamb$ be a maximal order in $D$. 

\noindent $1$. For any integral element $a\in \mc{O}_F^{\times}$ that is positive at all real archimedean places ramified in $D$, there exists an integral element $\al\in D^{\times}$ such that $\Nr(\al) = a$.

\noindent $2$. For any two left ideals $I$, $J$ of $\Lamb$ such that $\Nr(I) = a\Nr(J)$ for some element $a\in F^{\times}$ that is positive at all real archimedean places that $D$ is ramified at, there exists $\al\in D^{\times}$ such that $I\al = J$.\end{Thm}

Furthermore, it is well-known that for any division algebra $D$ central over a number field $F$, there exist infinitely many maximal subfields $L$ such that $L/F$ is a Galois extension. In fact, by the Grunwald-Wang theorem, there exists a maximal subfield $L\sub D$ such that $L/F$ is cyclic. Hence, one may ask whether the analog of the Hasse-Mass-Schilling theorem restricted to elements contained in Galois extensions holds.

\noindent \textbf{Question:} For what elements $a\in F^{\times}$ do there exist elements $\al\in D^{\times}$ such that the Galois closure of $F(\al)$ over $F$ has an embedding in $D$ and $\Nr(\al) = a$?

When the division algebra $D$ is the endomorphism algebra of a simple abelian variety, the corresponding endomorphisms of the abelian variety will be Galois isogenies, which is our primary interest exploring this. For quaternion algebras that fulfill Eichler's condition (i.e. not totally definite), it is clear that this is the case since all quadratic extensions are Galois. On the other hand, it is easy to see that this does not hold for all division algebras, as evidenced by the following counter-example. 

\begin{Example} \normalfont Choose distinct primes $p$, $l$ such that $l\not\equiv 1\pmod{p}$. Let $D$ be any central division algebra of dimension $p^2$ over $\bq$ with Hasse invariant $1/p$ at the prime $l$. Then $D$ fulfills Eichler's condition and $l$ is positive at the archimedean place of $\bq$. We show that $l$ does not occur as the reduced norm of any element that lies in a Galois extension of $\bq$ embedded in $D$.

Note that since $D$ has center $\bq$, for any element $\al\in D^{\times}$, we have $\Nr(\al) = \Nm_{\bq(\al)/\bq}(\al)$. Suppose there exists a Galois extension $L/\bq$ embedded in $D$ such that $l = \Nr(\al)$ for some $\al\in L^{\times}$. Since the extension $L/\bq$ is Galois of prime degree, it is cyclic. Furthermore, $L$ is a maximal subfield of $D$ and hence, $l\in \Nm_{L/\bq}(\mc{O}_L^{\times})$. Since $L/\bq$ has local degree $p$ at $l$, it follows that $l$ is totally ramified or inert in $L/\bq$. In the latter case, $l\mc{O}_L$ could not be a norm of an ideal of $L$. So $l$ is totally ramified in $L/\bq$ and the prime of $L$ lying above $l$ is principal. From the Kronecker-Weber theorem, it follows that $L\sub \bq(\zeta_{_N})$ for some integer $N$. Let $n$ be the smallest integer such that $L\sub \bq(\zeta_n)$. Then each prime dividing $n$ is either $p$ or is $\equiv 1 \pmod{p}$ and in particular, $l\nmid n$. Thus, $l$ is unramified in $\bq(\zeta_n)$ and hence, is unramified in $L$, a contradiction.\end{Example}

\begin{subsection}{\fontsize{11}{11}\selectfont Decomposition of the torsion subgroup schemes}\end{subsection}

Let $B$ be an abelian variety with endomorphism algebra $D$ of type $\IV(1,d)$ over a finite field $\bFq$. For any integer $N$ prime to $q$, the $N$-torsion subgroup scheme $B[N]$ is an etale $\Gal_{\bFq}$-module and $B(\bFP)[N]\cong (\bz/N\bz)^{2d}$. In this subsection, we explore the decomposition of this subgroup scheme into those annihilated by certain isogenies. We will need the following lemma (generally attributed to Iwasawa), a proof of which may be found in ([Was82], Chapter 10).

\begin{Lem} Let $p$ be a prime and let $L/K$ be a Galois extension of number fields with $[L:K]$ a $p$-power. If $L$ is ramified at precisely one finite prime, then $p$ divides the class number $\mr{h}(L)$ if and only if it divides $\mr{h}(K)$ .\end{Lem}

For a division algebra $D$ of dimension $n^2$ over a number field $F$, an extension $L/F$ of degree $n$ has an embedding in $D$ if and only if it fulfills the following equivalent conditions:

\noindent - $L$ splits $D$.\\
- For any place $\wti{v}$ of $L$ and the place $v$ of $F$ lying under $\wti{v}$, $L_{\wti{v}}$ splits the simple algebra $D_v:= D\otimes_{F} F_v$.

For any such field $L$, its ring of integers $\mc{O}_L$ is a finitely generated $\mc{O}_K$-submodule of $D$ and by ([AG60], Proposition 1.1) has an embedding in \textit{some} maximal order $\Lamb$ of $D$. Consequently, $\mc{O}_L$ has an embedding in every maximal order of $D$ conjugate to $\Lamb$. However, the question of which maximal orders of $D$ contain a copy of $\mc{O}_L$ is substantially more subtle and as far as we know, remains open in its full generality. A partial answer to this is provided by the following elegant theorem of Arenas-Carmona.

\begin{Thm} $\mr{([Car03])}$ Let $D$ be a central division algebra of dimension $n^2$ over a number field $F$ such that:

\noindent - $n\geq 3$.\\
- At any place $\mfp$ of $F$ where $D$ is ramified, the central simple algebra $D_v:= D\otimes_F F_v$ is a division algebra.

Then there exists an abelian extension $\Sigma_F$ of $F$ with the following property: for any degree $n$ extension $L$ of $K$ that splits $D$, the ring of integers $\mc{O}_L$ has an embedding in precisely $[\Sigma_F\cap L:F ]^{-1}$ of the conjugacy classes of maximal orders of $D$.\end{Thm}

The field $\Sigma_F$ is a spinor class field of $F$. For details, we refer the reader to [Car03]. In particular, if $L$ is a degree $n$ extension with an embedding in $D$ and the field extensions $L,\Sigma_F$ are linearly disjoint over $F$, then the ring of integers $\mc{O}_L$ has an embedding in every maximal order of $D$. On the other hand, if $L\subseteq \Sigma_F$, then $L$ is contained in precisely $\frac{1}{n}$-th of the conjugacy classes of the maximal orders.

The second condition on $D$ is equivalent to the Hasse invariant at each ramified place being of the form $j/n\in \bq/\bz$ for some integer $j$ relatively prime to $n$. Equivalently, for any maximal subfield $L\sub D$ and a prime $v$ of $F$ such that $D$ is ramified at $v$, there is a unique prime $\wti{v}$ of $L$ lying over $v$. Clearly, any abelian variety $B$ of type $\IV(1,d)$ has an endomorphism algebra that fulfills both of the conditions. Hence, the last theorem is directly applicable to our setting. We will need to impose the following conditions for some of the results in this section.




\begin{Def} \normalfont We say a division algebra $D$ of type $\IV(1,d)$ fulfills the condition $(***)$ if:

\noindent $\bullet$ $d$ is a prime power.\\
$\bullet$ The center $K$ does not lie in $\bq(\zeta_d)$.\\
$\bullet$ The class number of the center $K$ of $D$ is prime to $d$.\end{Def}

\begin{Def} \normalfont We say an abelian variety of type $\IV(1,d)$ fulfills the condition $(***)$ if its endomorphism algebra fulfills the condition $(***)$.\end{Def}

\noindent The condition that $d$ is a prime power is, by far, the most restrictive of the three. Since $\bq(\zeta_d)$ has a unique quadratic subfield, the assumption that $K\nsubseteq \bq(\zeta_d)$ is mild. We provide the following two examples to demonstrate that the third condition is fairly generic as well.

\begin{Example} \normalfont Consider the case where $d$ is a $2$-power larger than $2$. It is well-known that the class number of an imaginary quadratic field $K = \bq(\sqrt{-a})$ is odd if and only if $a$ is a prime $\equiv -1\pmod{4}$. Choose any such imaginary quadratic field $K$ and choose a rational prime $p$ that splits in $K$. Let $\mfp,\ov{\mfp}$ be the primes of $K$ lying over $\bq$ and let $h$ be the smallest integer such that $\mfp,\ov{\mfp}$ are principal. Let $D$ be the division algebra central over $K$ with Hasse invariants \vspace{-0.1cm}$$\inv_{\mfp}(D) = \frac{1}{d},\;\;\inv_{\ov{\mfp}}(D) = \frac{d-1}{d}.$$ Then we may construct an abelian variety over the field $\bF _{p^h}$ with endomorphism algebra $D$.\end{Example}

\begin{Example} \normalfont For the case when $d$ is a $p$-power for some odd prime $p$, it was shown by Hartung ([Har74]) that there exist infinitely many imaginary quadratic fields with class number not divisible by $p$. His proof makes use of the Kronecker relation: \vspace{-0.1cm}$$\mr{h}(4n-s^2) = 2\sum\limits_{r|n,\;r>\sqrt{n} } r $$ where $\mr{h}(N)$ denotes the class number of quadratic forms of discriminant $-N$. Furthermore, by ([KO99], Theorem 1.1), the number of imaginary quadratic fields with discriminant $\leq X$ and class number indivisible by $p$ is bounded below as follows: 

For any constant $\epsilon > 0$, there exists a constant $X_{\epsilon}$ such that \vspace{-0.1cm}$$\#\{ -X_{\epsilon} < D <0:\; p\nmid \mr{h}(D)  \}  \geq \left(\frac{2(p-2)}{\sqrt{3}(p-1)}-\epsilon\right) \frac{\sqrt{X_{\epsilon}}}{\log(X_{\epsilon})},$$ which implies that division algebras fulfilling condition $(***)$ are ubiquitous. As seen in Proposition 3.7, for any such division algebra $D$, we may construct an abelian variety whose endomorphism algebra is isomorphic to $D$.
\end{Example}

\begin{Thm}  Let $B$ be an abelian variety of type $\IV(1,d)$ over a finite field $\bFq$ with the endomorphism ring $\End(B)$ a maximal order in the endomorphism algebra $D:=\End^0(B)$. Let $n$ be an integer relatively prime to $q$ that lies in the image of the norm map $\Nm_{K/\bq}: K\lra \bq$. 

\noindent $(1)$. There exists a Galois quasi-isogeny $\psi:B\lra B$ of degree $n$. 

\noindent $(2)$. Furthermore, if $B$ fulfills condition $(***)$ and $\End(B)$ is a maximal order in $D$, then for any sufficiently large integer $e$, there exist distinct, pairwise commutative Galois endomorphisms $\psi_i:B\lra B,\; i=1,\cdots, d$ of degree $n^e$.\end{Thm}


\begin{prf} $(1).$ Let $l_1,\cdots,l_r$ be the distinct primes dividing $n$. We show that there pairwise commutative exist Galois quasi-isogenies $\phi_1,\cdots,\phi_r$ with $\Nr(\phi_i) = l$. Note that pairwise commutativity is equivalent to the quasi-isogenies all lying in some maximal subfield of $D$.

We choose a prime $v$ of $\bq(\zeta_d)$ satisfying all of the following conditions. 

\noindent (i). $v$ has local degree one over $\bq$.\\
(ii). $v$ splits completely in $\bq(\zeta_{d},\sqrt[d]{l}_1,\cdots,\sqrt[d]{l}_r)$.\\
(iii). $v$ is inert in the degree $d$ cyclic extensions $\bq(\zeta_{d^2})$ and in $\bq(\zeta_{d},\sqrt[d]{p})$.\\
(iv). $v$ is prime to the discriminant of $K$.\\
(v). $v$ is inert in the quadratic extension $K(\zeta_d)$.

Note that the existence of such a prime $v$ is guaranteed by the Chebotarev density theorem. In fact, the set of the primes of $\bq(\zeta_d)$ fulfilling all of these conditions has Dirichlet density \vspace{-0.1cm}$$\frac{\phi(d)^2}{2d^{r+2}} > 0.$$ Let $L_0$ be the unique subfield of $\bq(\zeta_{{p_{_1}}})$ of degree $d$. Set $L:=L_0K$, so that $L$ is a CM field of degree $2d$ cyclic of degree $d$ over $K$.

Let $p_{_1}$ be the rational prime lying under $v$. By condition (iii), the inertia degree of $p$ in $\bq(\zeta_{{p_{_1}}})/\bq$ divisible by $d$ and $\gcd(d,\frac{p-1}{d}) = 1$. Hence, $p$ is inert in the degree $d$ subfield $L_0/\bq$ of $\bq(\zeta_{{p_{_1}}})$. Furthermore, since $p$ splits (completely) in $K/\bq$ and is inert in $L_0/K$, it follows that the two primes of $K$ lying over $p$ are inert in $L/K$. So, in particular, $L/K$ has local degree $d$ at the primes of $K$ that $D$ is ramified at and hence, $L$ splits $D$. Thus, $L$ has an embedding in $D$.

By condition (ii), the primes $l_1,\cdots, l_r$ split completely in $L_0$. Since these primes split completely in $K$, it follows that they split completely in the compositum $L = L_0K$. It suffices to show that each $l_i$ occurs as the degree of a quasi-isogeny. Hence, we fix a prime divisor $l$ of $n$ for the rest of this proof. 

Furthermore, since there exist infinitely many primes $v$ fulfilling the conditions in the last paragraph, we may choose $v$ such that $p_{_1}:=\charac(v)$ is unramified in the spinor class field $\Sigma_K$. By Arenas-Carmona's theorem, it follows that the ring of integers $\mc{O}_L$ has an embedding in every maximal order of $D$. In particular, $\mc{O}_L$ has an embedding in $\End(B)$.

Now, if any element $\al\in L$ has norm $\Nm_{L/K}(\al) = l$, then $\al$ corresponds to a quasi-isogeny $\psi_{\al}:B\lra B$ with degree $l$. We first show that $l$ is a local norm at all places of $K$ in the extension $L/K$. Let $v$ be any finite place of $L$ and let $w$ be the place of $K$ lying under $v$. We will need to show that $l$ lies in the image of the norm map $$\Nm_{L_v/K_w}:L_v\lra K_w $$ of local fields. We treat the following three cases separately.

\noindent \underline{Case 1} $\charac(v) = l$.

Since the primes of $K$ lying over $l$ split completely in $L$, we have $L_v \cong K_w$ and hence, $l$ is trivially in the image of the norm map.

\noindent \underline{Case 2} $\charac(v) = p_1$.

Since $L/K$ is cyclic and totally ramified at the primes of $K$ lying over $p_1$, it follows that $L_v/K_w^{\times}$ is cyclic of degree $d$. By construction, $l$ is a $d$-th power in $\bF _{p_1}^{\times}$. Hence, by Hensel's lemma, there exists $\al\in K_w^{\times}$ such that $l= \al^d$. Thus, $l = \al^d = \Nm_{L_v/K_w}(\al)$.

\noindent \underline{Case 3} $\charac(v) \notin \{l,p_1\}$.

Since $L_v/K_w$ is unramified, the norm map is surjective on the group of units in $K_w^{\times}$. Since $l$ is a unit in $K_w^{\times}$, there exists some $u\in L_v^{\times}$ such that $\Nm_{L_v/K_w}(u) = l$.

Thus, $l$ is locally a norm at all places of $L/K$ and since $L/K$ is cyclic, it follows from the Hasse norm theorem that $l$ lies in $\Nm_{L/K}(L^{\times})$. Let $\al\in L^{\times}$ be an element such that $\Nm_{L/K}(\al) = l$. Then any quasi-isogeny corresponding to $\al$ has degree $l$.

\noindent $(2).$ Now suppose, $B$ fulfills the condition $(***)$. We construct the number field $L$ in the same manner as part $(1)$. Now, $L/K$ is cyclic of degree $d$ which is assumed to be a prime power. Furthermore, the rational prime $p_1$ is inert in $K/\bq$. So $L/K$ is totally ramified at the prime $p_1\mc{O}_K$ and is unramified everywhere else. Since the class number of $K$ is relatively prime to $d$, by Lemma 3.26, so is the class number of $L$. Let $\wti{\mfl}$ be a prime of $L$ lying over $\mfl$ and let $h$ be the smallest integer such that $\wti{\mfl}^h$ is principal, with $\wti{\mfl}^h = \al\mc{O}_L$ for some $\al\in\mc{O}_L$. Now, $h$ divides the class number of $L$ and hence, is relatively prime to $d$. Hence, for any integer $N \geq dh - d- h$, there exist non-negative integers $a_1, a_2$ such that $a_1h+a_2d = N$. Thus, we have $$\mfl^N = \Nm_{L/K} (\al^{a_1}\mfl^{a_2}) = \Nr(\al^{a_1}\mfl^{a_2}),$$ which completes the proof.\end{prf}


\noindent \textbf{Remark} From the proof, it is clear that we may construct infinitely many such endomorphisms of degree $l^N$, no two of which commute. Any two such endomorphisms generate $D$ as a $\bq$-algebra.

\begin{Corr} For any sufficiently large integer $N$, there exist cyclic endomorphisms $\psi_1,\cdots,\psi_d: B\lra B$ and an integer $m<N$ such that:

\noindent - the $\psi_i$ are pairwise commutative.\\
- the endomorphisms $l^m\psi_j$ are of degree $l^N$.\\ 
- $B[l^N] = \sum\limits_{i=1}^d B[\psi_i]$.\\
- $B[\psi_j] \cap B[\what{\psi}_j] = \{0_B\}$ where $\what{\psi}_j = \prod\limits_{i\neq j}\psi_i$.\end{Corr}

\begin{prf} As shown in the preceding theorem, there exist infinitely many cyclic endomorphisms $\phi:B\lra B$ of degree $l^N$ such that $K(\psi)/K$ is a Galois extension. Let $\psi$ be one such endomorphism and let $\psi_1,\cdots,\psi_d$ be the Galois conjugates of $\psi$. Since $K(\psi)/K$ is Galois, the $\psi_j$ lie in the embedding of $K(\psi)$ in $D$ . Hence, they are pairwise commutative.

Furthermore, the ideal $(\psi_i,\psi_j) = (1)$ in $\mc{O}_{K(\psi)}$. Hence, $\ker(\psi_i)\cap \ker(\psi_j) = \{0 \} $.\end{prf}

\begin{Prop} Let $B_1$, $B_2$ be simple abelian varieties over a finite field $\bFq$ with the same endomorphism algebra $D$ of type $\IV(1,d)$ and $\Lamb_i:=\End(B_i)$ maximal orders in $D$ for $i=1,2$. There exists a cyclic isogeny $\psi:B_1\lra B_2$.\end{Prop}

\begin{prf} Since $B_1$ and $B_2$ are of type $\IV(1,d)$ and have the same endomorphism algebra, they are isogenous by Proposition 2.2. Now, there exists a left ideal $I$ of $\Lamb_1$ such that $\phi_I(B_1) = B_2$. Choose a prime ideal $\mfp$ of $K$ prime to $N$ such that: 

\noindent - $\mfp$ has local degree one over $\bq$.\\
\noindent - $\mfp$ lies in the same ideal class as $\Nr(I)$ in $\mr{Cl}(K)$. 

Choose a maximal left ideal $J$ of $\End(B_1)$ containing $\mfp$. By Eichler's theorem, $J$ is in the same ideal class as $I$ in the ideal class set $\Cl(\Lamb_1)$. Hence, $\phi_{J} = \phi_{I}$ up to isomorphism.\end{prf}



\begin{Prop} Let $B_1,B_2$ be abelian varieties with endomorphism algebra $D$ of type $\IV(1,d)$ with the endomorphism rings $\Lamb_1:=\End(B_1)$, $\Lamb_2:=\End(B_2)$ maximal orders in $D$. The following are equivalent.

\noindent $1.$ There exist a prime ideal $I$ of $\mc{O}_K$ such that $\phi_I(B)$, $B_2$ are Galois conjugates.\\
$2$. $\Lamb_2 = \al^{-1}\Lamb_1 \al$ for some $\al\in D$.\end{Prop}

\begin{prf} $(1)\Leftarrow (2)$: Suppose $B_2 = \sigma(B_1)$ for some $\sigma\in \Gal(\bFq/\bFp)$. Clearly, $\sigma\circ I_0$ yields an isomorphism $\Lamb_1\cong \Lamb_2$ which, when tensored with $K$, yields a automorphism of the division algebra $D$. By the Skolem-Noether theorem, this is an inner automorphism by some element of $D$. Hence, $\Lamb_1$, $\Lamb_2$ are conjugate maximal orders. 

\noindent $(2)\Rightarrow (1)$: Conversely, suppose $\Lamb_2 = \al^{-1}\Lamb_2 \al$ for some $\al\in D^{\times}$. Now, there exists a left ideal $I$ of $\Lamb_1$ such that $\phi_I(B_1) = B_2$ and $O_r(I) = 	\Lamb_2$. Hence, $O_r(I\al) = \al^{-1}O_r(I)\al$, $O_l(I\al) = O_l(I)$. Replacing $B_2$ by another abelian variety in its isomorphism class if necessary, we may assume $I$ is an integral two sided ideal of $\Lamb_1$. Furthermore, $\phi_I$ is a composition \vspace{-0.1cm}$$\phi_{\mr{sep}} \circ\phi_{\mr{insep}}: B_1\lra B_3\lra B_2$$ of a separable and a purely inseparable morphism. By ([Yu11], theorem 1.3), there exists an abelian variety $\wti{B}_3$ such that $\End(\wti{B}_3)$ is a maximal order, $End(B)\sub \End(\wti{B}_3)$ and every isogeny to $B_3$ factors through $\wti{B}_3$. Hence, there is a composition \vspace{-0.1cm}$$\wti{\phi}_{\mr{sep}} \circ\wti{\phi}_{\mr{insep}}: B_1\lra \wti{B}_3\lra B_2$$ where $\wti{\phi}_{\mr{sep}}$ is purely inseparable and $\wti{\phi}_{\mr{insep}}$ and $\wti{\phi}_{\mr{sep}}$ is separable.  

Since the abelian varieties $\wti{B}_3, B_2$ have endomorphism rings that are maximal orders, the isogeny $\wti{\phi}_{\mr{sep}}$ linking them is an ideal isogeny, i.e. there exists a maximal left ideal $J$ of $\End(\wti{B}_3)$ such that $\wti{\phi}_{\mr{sep}} = \phi_J$ up to isomorphism. Since $\phi_J$ is separable, $\Nr(J)$ is relatively prime to $p$.

Since $\phi_{\mr{insep}}$ is some power of the Frobenius, it suffices to show that the separable morphism $\phi_{\mr{sep}}$ is an isomorphism. Thus, we may assume $\phi_I$ is separable, meaning $\Nm(I)$ is prime to $p$. Since $D$ is split at every prime of $K$ dividing $I\cap \mc{O}_K$, it follows that $I = I_0\Lamb_1$ for some ideal $I_0$ in $\mc{O}_K$. Hence, $\phi_{\mr{sep}}$ is an isomorphism and $B_2 = \phi_{\mr{insep}}(B_1)$.\end{prf}

\noindent In the case of supersingular elliptic curves, it is well-known that the endomorphism ring is a maximal order in the quaternion algebra $\bq_{p,\infty}$ ramified only at $p$ and the archimedean place of $\bq$. In our setting of an abelian variety of type $\IV(1,d)$, we do not quite have the analogous result since the center \vspace{-0.1cm}$$\Cent(\End(B)) = \End(B)\cap \Cent(D)$$ of $\End(B)$ might not be the maximal order $\mc{O}_K$. However, we show that the failure of $\End(B)\cap K$ to be the maximal order in $K$ is a good measure of the failure of $\End(B)$ to be a maximal order in $D$. Note that in general, the set of orders of $D$ is rather complicated and this does not easily generalize to arbitrary abelian varieties over finite fields.  

\begin{Prop} Let $B$ be an abelian variety of type $\IV(1,d)$ over a finite field of characteristic $p$. Let $K$ be the center of $D:=\End^0(B)$ and let $\mc{O}$ be the center of $\End(B)$.

\noindent $1$. The center of $\End(B)$ is isomorphic to an order $\mc{O}$ in the field $K$.

\noindent $2$. The order $\Lamb:= \End(B) \otimes_{\mc{O}} \mc{O}_K$ is a maximal order in $D$. In particular, $\End(B)$ is a maximal order in $D$ if and only if $\End(B)\cap K = \mc{O}_K$. \end{Prop}

\begin{prf} 1. We have $\Cent(\End(B))\otimes_{\bz} \bq\cong D$ and since the centralizer $\Cent(\End(B))$ is a finitely generated $\bz$-module, it follows that it is an order in $K$.

\noindent 2. By ([AG60], Proposition 1.2), it suffices to show that the localization of $\Lamb$ at any prime $v$ of $K$ is a maximal order in $D\otimes_{K} K_v$. Let $\mfl$ be a prime of $K$ that does not divide $p$. By Tate's isogeny theorem, we have an isomorphism \vspace{-0.1cm}$$\End(B)\otimes_{\mc{O}_K} \mc{O}_{K_{\mfl}} \cong \Mat_d(\mc{O}_{K_{\mfl}})$$ which is a maximal order in $\Mat_d(K_{\mfl})$.  


Let $v$ be a prime of $K$ lying over $p$ and let $\frac{j}{d}$ be the Hasse invariant of $D$ at $v$, for some integer $j$ prime to $d$. Since $p$ splits completely in $K$, we have $K_v\cong \bq_p$. The Newton polygon of $B$ is $d\times \frac{j}{d}, d\times \frac{d-j}{d}$ for some integer $j$ prime to $d$. Now, $D_v:=D\otimes_{K} K_v$ is a division algebra of degree $d^2$ over $\bq_p$ with Hasse invariant $j/d$. Furthermore, by Tate's isogeny theorem, the subring $\Lamb\otimes_{\mc{O}_K} \mc{O}_{K_v}$ of $D_v$ is the endomorphism ring of the pure isocrystal of slope $j/d$. By a theorem of Dieudonne-Lubin ([Fro68], Page 72), this is the (unique) maximal order in the division algebra over $\bq_p$ with Hasse invariant $j/d\in \bq/\bz$.

Thus, we have verified that the localization of $\Lamb$ at any prime is a maximal order. Hence, $\Lamb$ is a maximal order in $D$.\end{prf}

\begin{Corr} Let $B$ be an abelian variety of type $\IV(1,d)$ over a finite field and let $\pi$ be the Weil number associated to it. Let $\Lamb$ be a maximal order of $D:=\End^0(B)$ that contains $\End(B)$. Then the index $[\Lamb:\End(B)]$ divides the index $[\mc{O}_K: \bz[\pi,\ov{\pi}]]$ where $K$ is the center of $D$.\end{Corr}

\begin{prf} Let $\mc{O}$ be the center of $\End(B)$. As shown in the preceding proposition, the index $[\Lamb:\End(B)]$ coincides with the index $[\mc{O}_K:\mc{O}]$. Now $\mc{O}$ contains the order $\bz[\pi,\ov{\pi}]$ and hence, $[\mc{O}_K:\mc{O}]$ divides $[\mc{O}_K: \bz[\pi,\ov{\pi}]]$.\end{prf}

\begin{Corr} Let $B_1$, $B_2$ be abelian varieties over a finite field with endomorphism algebra $D$ of type $\IV(1,d)$. Let $\Lamb_1$, $\Lamb_2$ be maximal orders of $D$ containing the orders $\End(B_1)$, $\End(B_2)$ respectively. Then there exists an ideal isogeny linking $B_1$ and $B_2$ if and only if \vspace{-0.1cm}$$[\Lamb_1:\End(B_1)] = [\Lamb_2:\End(B_2)].$$ \end{Corr}

\begin{Prop} Let $D$ be a division algebra of type $\IV(1,d)$ with center $K$ and let $B, B'$ be abelian varieties over a finite field $\bFpk$ with endomorphism algebra $D$. Let $\phi:B\lra B'$ be a separable isogeny of a prime degree $l$. Let $\Lamb := \End(B)$, $\Lamb' := \End(B')$ be the endomorphism rings with centers $\mc{O}$, $\mc{O}'$ respectively. Then the following are equivalent:

\noindent $1$. $\mc{O} = \mc{O}'$.

\noindent $2$. The left ideal $I(\ker(\phi)):=\{\psi\in \mc{O}:\;\psi(\ker(\phi)) = 0 \}$ is a locally free module of $\mc{O}$ of norm equal to $\deg(\phi)$.

\noindent $3$. There exists an isogeny $B\lra B'$ of degree prime to $\deg(\phi)$.\end{Prop}

\begin{prf} ($1\Leftrightarrow 2$): Let $*$ denote the involution of the second kind that $D$ is equipped with. Then \vspace{-0.1cm}$$\bz + l^2\mc{O}\sub \bz+\phi\mc{O}'\phi^* \sub \mc{O}.$$ By Tate's isogeny theorem, we have \vspace{-0.1cm}$$\Hom(B,B')\otimes_{\bz}\bz_l \cong \Hom_{\bz[\Fr_k]}(T_l(B),T_l(B')).$$ Hence, the Tate modules $T_l(B)$, $T_l(B')$ are isomorphic as $\bz_l[\Fr_k]$-modules if and only if the orders $\Lamb_l$, $\Lamb'_l$ are isomorphic, which in turn is equivalent to the centers $\mc{O}_l$, $\mc{O}'_l$ being isomorphic. Since \vspace{-0.1cm}$$I(\ker(\phi))\otimes_{\mc{O}}\mc{O}_{\mfp} = \mc{O}_{\mfp}$$ for any prime $\mfp$ of $\mc{O}$ lying over $p$, this is equivalent to $I(\ker(\phi))$ being locally free as a left $\mc{O}$-module. 

\noindent $(2\Rightarrow 3)$: If $I(\ker(\phi))$ is a locally free $\mc{O}$-module, we may choose an ideal $J$ of $\mc{O}$ prime to $I(\ker(\phi))$ and in the same ideal class as $I(\ker(\phi))$. The ideal isogeny $\phi_J: B\lra B'$ is of degree prime to $\deg \phi$. 

\noindent $(3\Rightarrow 1)$: Since one of $\mc{O}$ and $\mc{O}'$ is contained in the other, we may assume without loss of generality that $\mc{O}'\sub \mc{O}'$. So $[\mc{O}:\mc{O}']$ divides $l$. Let $\phi_1$ be an isogeny of degree prime to $l$. Then $[\mc{O}:\mc{O}']$ divides $(\deg\; \phi_1)^d$ and since $[\mc{O}:\mc{O}']$ also divides $l$, it follows that $\mc{O} = \mc{O}'$.\end{prf}

\noindent The following corollary is immediate.

\begin{Corr} Let $D$ be a division algebra of type $\IV(1,d)$ with center $K$ and let $B, B'$ be abelian varieties over a finite field $\bFpk$ with endomorphism algebra $D$. Let $\phi:B\lra B'$ be a separable isogeny of degree a prime $l$ . Let $\Lamb, \Lamb'$ be maximal orders of $D $ containing the endomorphism rings $\End(B)$, $\End(B')$. Then the following are equivalent:

\noindent $1$. $[\Lamb:\End(B)] = [\Lamb':\End(B')]$.

\noindent $3$. There exists an isogeny $B\lra B'$ of degree prime to $\deg(\phi)$.\end{Corr}

\noindent Earlier, we showed that an isogeny $B\lra B'$ of abelian varieties of type $\IV(1,d)$ arising from a left ideal of $\End(B)$ is a normed isogeny. We now show that the converse is also true.

\begin{Prop} Let $B$ be an abelian variety of type $\IV(1,d)$ with $\End(B)$ a maximal order in $D:=\End^0(B)$. Any ideal isogeny from $B$ to another abelian variety is given by $\phi_{I}$ for some left ideal $I$ of $\End(B)$.\end{Prop}

\begin{prf} Write $\Lamb:=\End(B)$. As seen in the last proposition, $\Lamb$ is a maximal order in $D$. It suffices to show that for any ideal $n$ of $K$, the number of subgroup schemes of $B$ order $n$ coincides with the number of left integral ideals of reduced norm $n$. Furthermore, it suffices to show this for the case when $n$ is a prime power.

\noindent \textbf{Case 1.} $n = p^m$ for some integer $m$.

Since $B$ has $p$-rank zero, every isogeny $\phi:B\lra B'$ of order $p^m$ is purely inseparable and factors through $B^{(p^m)} := \Fr_B^{p^m}(B)$. Since the degrees coincide, we have $B'\cong B^{(p^m)}$ and there is a unique subgroup-scheme of order $p^m$ in $B[p^m]$. On the other hand, $\Lamb_{\mfp}$ is the unique maximal order in a $d^2$-dimensional central division algebra over $\bq_p$ and hence, has a unique ideal of reduced norm $p^m$.

\noindent \textbf{Case 2.} $n = l^m$ for some integer $m$ and prime $l\neq p$.

Let $\mfl$ be a prime of $l$ lying over $l$. We have $\Lamb_{\mfl}\cong \Mat_d(\mc{O}_{K_{\mfl}})$ and hence, the left ideals of $\Lamb_{\mfl}:= \Lamb\otimes_{\mc{O}_{K}} \mc{O}_{K_{\mfl}}$ are principal. Now, the set of elements of $\Mat_d(\mc{O}_{K_{\mfl}})/\GL_d(\mc{O}_{K_{\mfl}})$ of determinant $l^m$ is in bijection with the set of submodules of $\mc{O}_{K_{\mfl}}^d$ of index $l^m$. For any such submodule $M\sub \mc{O}_{K_{\mfl}}^d$, the image of $M$ under the homomorphism: $$ 0\lra l^m\mc{O}_{K_{\mfl}}\lra \mc{O}_{K_{\mfl}}^m\lra (\bz/l^m\bz)^d $$ is a subgroup of order $l^m$. Conversely, let $H$ be a subgroup of order $l^m$ in $(\mc{O}_{K}/\mfl^m\mc{O}_{K})^d$. The preimage of $H$ in $\mc{O}_{K_{\mfl}}^d$ is a $\mc{O}_{K_{\mfl}}$-submodule of $\mc{O}_{K_{\mfl}}^d$ with a cokernel of order $l^m$.\end{prf}

\begin{Prop} Let $B$ be an abelian variety of type $\IV(1,d)$ with $\End(B)$ a maximal order in $\End^0(B)$. Any normed isogeny from $B$ to another abelian variety is given by $\phi_{J}$ for some maximal left ideal $J$ of $\End(B)$ with the degree $\deg \phi_I$ a prime.\end{Prop}

\begin{prf} As shown in the last proposition, the isogeny is of the form $\phi_{I}$ for some left ideal $I$ of $\End(B)$. Let $\Nr(I)\sub \mc{O}_K$ be the reduced norm ideal of $I$. Let $\mr{H}(K)$ denote the Hilbert class field of $K$. Applying the Chebotarev density theorem to the abelian extension $\mr{H}(K)/K$, we see that in every ideal class of $K$, there exist infinitely many primes of $K$ with local degree one over $\bq$. Choose any such prime $\mf{q}$ in the same ideal class as $\Nr(I)$ and let $J$ be a maximal left ideal of $\End(B)$ containing $\mf{q}$. By Eichler's theorem, $J$ lies in the same ideal class as $I$. Thus, $\phi_I = \phi_J$ up to isomorphism.\end{prf}

\begin{Lem} Let $D$ be a central division algebra over a global field $F$, $\Lamb$ a maximal order and $I$ a left ideal in $\Lamb$. If $J \sub I$ is another left $\Lamb$-ideal and $\Nr(J) = \Nr(I)$, then $J = I$. In particular, if there exists an element $\al \in I$ such that $\Nr(\al) = \Nr(I)$, then $I = \Lamb \al$.\end{Lem}
\begin{prf} See ([Rei75], Chapter 5).\end{prf}






\begin{Prop} Let $B$ be an abelian variety of type $\IV(1,d)$ fulfilling the condition $(***)$. Write $D:=\End^0(B)$ and let $K$ be its center. Let $l$ be a rational prime that splits completely in the Hilbert class field $\mr{H}(K)$. For any sufficiently large integer $N$, there exist endomorphisms $\psi,\tau\in \End(B)$ and an integer $m$ such that:

\noindent $1.$ $\psi$ is a cyclic isogeny of degree $l^{N-dm}$.\\
$2.$ The isogenies $\psi_j:= \tau^{-j}\psi\tau^j$ $(j=0,\cdots,d-1)$ are pairwise commutative and yield the decompositions\vspace{-0.1cm}$$B[l^{N-dm}] = \bigoplus\limits_{j=0}^{d-1} B[\psi_j],\;\;\;B[l^{N}] = \sum\limits_{j=0}^{d-1} B[l^m\psi_j]$$ of the torsion subgroup schemes.\\
$3.$ $K(\psi, \tau)  = D$, i.e. $\psi, \tau$ generate $D$ as a $K$-algebra..  \end{Prop}

\begin{prf} 
Since $\mc{O}_L$ is a finitely generated $\mc{O}_K$-module, Proposition 1.1 of [AG60] implies that it is contained in some maximal order $\Lamb$ of $D$. Let $\sigma$ be a generator of the Galois group $\Gal(L/K)$. By the Skolem-Noether theorem, $\sigma$ extends to an inner automorphism of $D$ by an element $\al\in D^{\times}$. Replacing $\al$ by an integer multiple if necessary, we may assume without loss of generality that $\al\in \Lamb^{\times}$. 

Furthermore, $\tau^d$ commutes with $\psi$ and hence, $\tau^d$ lies in $L$. On the other hand, $\tau^j\notin L$ for $j= 1,\cdots,d-1$ and hence, the field $K(\tau)$ is of degree $d$ over $K$ and is disjoint from $L$ over $K$. Hence, the inclusion $K(\psi,\tau) \hra D$ is an equality.\end{prf}

\begin{Prop} Let $D$ be a division algebra of type $\IV(1,d)$ and let $\mc{O}_D$ be a fixed maximal order in $D$. There exists bijections: 

\noindent\begin{tikzcd} \{\text{Isomorphism classes of abelian varieties over } \bFP \text{ with endomorphism ring a maximal order in } D\} \arrow{d} \\ 
\{\text{Left ideal classes of } \mc{O}_D \ \} \arrow{d} \\ 
\{\text{Maximal left ideal classes of } \mc{O}_D \text{ of prime degree}\}. \end{tikzcd}\end{Prop}

\begin{prf} By Corollary 3.21, there exists an abelian variety with endomorphism ring $\mc{O}_D$. As seen above, the abelian varieties with endomorphism algebra $D$ define an isogeny class. Fix an abelian variety $B$ in this isogeny class and let $\Lamb$ be its endomorphism ring, a maximal order in $D$. Any isogeny $B\lra B_1$ is given by $\phi_I$ for some left ideal $I$. Furthermore, any two left ideals $I,J$ satisfy $\phi_I = \phi_J$ up to isomorphism if and only if they lie in the same ideal class. Since there is a bijection between the left ideal class sets of $\Lamb$ and $\mc{O}_D$, this yields the first bijection.

For the second bijection, let $I$ be any left ideal $\mc{O}_D$. Choose a prime $\mf{q}$ of $K$ in the same ideal class as $\Nr(I)$ in $\mr{Cl}(K)$. Let $J$ be a maximal left ideal of $\mc{O}_D$ containing $\mf{q}$. Then $\Nr(J) = \mf{q}$ and hence, the ideals $\Nr(I)$ and $\Nr(J)$ lie in the same ideal class of $K$. From Eichler's theorem, it follows that $I$ and $J$ lie in the same left ideal class of $\mc{O}_D$.\end{prf}

\begin{subsection}{\fontsize{11}{11}\selectfont The group of two-sided ideals}\end{subsection}

As observed a long time ago (see [Rei75]), the theory of two-sided ideals in orders of division algebras is vastly different from that of the one-sided ideals. In fact, while the one-sided ideals do not have a multiplicative group structure, it is easily verified that two-sided ideals of $\Lamb$ form an abelian group under multiplication. We say two-sided ideals $I$ and $J$ are \textit{equivalent} if $J= \al I$ for some $\al\in \Lamb^{\times}$. We denote the group of two sided ideals up to this equivalence by $\mr{Idl}(\Lamb)$. It is a well-known fact that is a finite abelian group. 

Note that a pair $I,J$ of two-sided ideals may be equivalent in the set of left ideal classes despite lying in different classes as two-sided ideals. However, the class group $\mr{Idl}(\Lamb)$ is not far from the class group of $K$, as shown in the next proposition.

\begin{Prop} Let $\Lamb$ be a maximal order in $D$. There exists a bijection \vspace{-0.1cm} $$\{\text{Two-sided prime ideals of } \Lamb \} \lra \{\text{Prime ideals of } \mc{O}_K \}, \;\;\; J\mapsto J\cap\mc{O}_K.\vspace{-0.1cm}$$ Furthermore, we have the exact sequence $$ 0\lra \mr{Cl}(\mc{O}_K)\lra \mr{Idl}(\Lamb)\lra (\bz/d\bz)^2\lra 0.\vspace{-0.1cm}$$
\end{Prop}

\begin{prf} Note that for any prime ideal $\mf{q}$ of $K$, $\Lamb \mf{q}\Lamb \cap K = \mf{q}$ and hence, the map is surjective. 

For the injectivity, let $\mfp$ be a prime ideal in $K$ and $P$ a two-sided prime ideal of $\Lamb$. Then the completion $P_{\mfp}$ is a maximal ideal of $\Lamb_{\mfp}$. If $\mfp$ is one of the two primes dividing $p$, then $\Lamb_{\mfp}$ is the maximal order in the a $d^2$-dimensional division algebra over $K_{\mfp}$. In this case, $\Lamb_{\mfp}$ has a unique two-sided maximal ideal $P_{\mfp}$ with $P_{\mfp}^d = \mfp\Lamb_{\mfp}$. On the other hand, if $\mfp\nmid p$, we have $\Lamb_{\mfp} \cong \Mat_d(\mc{O}_{K_{\mfp}})$. In this case, the only maximal two-sided ideal of $\Lamb_{\mfp}$ is $\mfp \Lamb_{\mfp}$. 

We have a natural homomorphism \vspace{-0.1cm}$$\mr{Idl}(\mc{O}_K)\lra \mr{Idl}(\Lamb),\; I\mapsto \Lamb I\Lamb.$$ To see that it is injective, suppose $\Lamb I\Lamb = \Lamb$ for some ideal $I$ of $\mc{O}_K$. Then $I^d = \Nr(I) = \mc{O}_K$ and hence, $I = \mc{O}_K$. The co-kernel is described by the last paragraph.\end{prf}

\begin{subsection}{\fontsize{11}{11}\selectfont A few examples of Jacobians of type $\IV(1,d)$}\end{subsection}

\noindent We know that any principally polarized abelian variety of dimension three occurs as the Jacobian of a smooth projective curve. Furthermore, this curve is either hyperelliptic or planar quartic. Unfortunately, for higher dimensions, we do not have a concrete characterization of such abelian varieties that arise as Jacobians of curves. In this subsection, we give a few examples of abelian varieties of type $\IV(1,d)$ that arise as Jacobians of curves.\\

\begin{Thm} $\mr{(Oort)}$ For any prime $p$ and integer $g\geq 3$, there exists a hyperelliptic curve over $\bFP$ whose Jacobian has $p$-rank zero.\end{Thm}

In particular, for any prime $p$, Oort's theeorem guarantees the existence of an abelian variety over $\bFP$ of type $\IV(1,3)$ that is the Jacobian of a hyperelliptic curve.



\begin{Example} \normalfont Let $C$ be the hyperelliptic curve $C: y^2 = 1 - x^7$ over $\bq$ and let $A:=\mr{Jac}(C)$ be its Jacobian. Then $A$ is a principally polarized abelian variety of dimension three over $\bq$ with endomorphism algebra $\End(A) = \bz(\zeta_7)$. Furthermore, $A$ has good reduction at all primes outside the set $\{2,7\}$ and since it has complex multiplication, it has potential good reduction everywhere. Let $p$ be a prime $\equiv 2\pmod{7}$. Then $p$ has inertia degree $3$ in $\bq(\zeta_7)$. By the theory of complex multiplication, the Newton polygon of the reduction $A_v$ is $3\times 1/3,\;3\times 2/3$. Hence, $\End^0(A_v)$ is a central division algebra $D$ of dimension $9$ over the imaginary quadratic field $\bq(\sqrt{-7})$ such that $\bq(\zeta_7)$ splits $D$. Furthermore, we have \vspace{-0.1cm}$$\bz(\zeta_7)\cong \End(A)\hra \End(A_v)$$ and hence, $\End(A_v)$ is a maximal order in $D$.\end{Example}

\begin{Example} \normalfont Let $D$ be a division algebra of type $\IV(1,3)$ with center $\bq(\sqrt{-1})$. Let $B$ be a principally polarized abelian variety over a finite field $\bFq$ with $\End^0(B) = D$. Then $B$ is the Jacobian of some curve $C$ over $\bFq$. Furthermore, since $B$ has a polarization-equivariant automorphism of order $4$, it follows from [Mil86] that the curve $C$ is hyperelliptic.\end{Example}

\begin{Example} \normalfont This construction using metacyclic Galois covers of the projective line draws from the results of [Ell01] and [CLS11]. Choose a rational prime $q\equiv 7\pmod{12}$. Set $d:=\frac{q-1}{6}$ and fix an integer $k$ of order $3$ modulo $q$. Consider the group \vspace{-0.1cm}$$G_{q,3}:= \la a,b:\; a^q = b^3 = 1,\; b^{-1}ab = a^k   \ra .$$ Let $Y$ be a Galois covering of the projective line $\mb{P}^1$ by the group $G_{q,3}$ over $\bc$ with three branched points $P_1,P_2, P_3$ each of ramification index $3$. Then the quotient curve $X:=Y/\la b\ra$ is a smooth projective curve and by Theorem 1 of [CLS11], its Jacobian $\Jac(X)$ is a simple abelian variety that admits complex multiplication by the unique index $3$ subfield $\bq(\zeta_q^{(3)})$ of $\bq(\zeta_q)$. In particular, $\Jac(X)$ is simple of dimension $\frac{q-1}{6}$. Since $\Jac(X)$ has complex multiplication, it follows from a theorem of Shimura that it has a model $A$ over some number field $F$. Enlarging $F$ if necessary, we may assume $F = F_{A}^{\conn}$ and in particular, $\bq(\zeta_q^{(3)})\sub F$. Choose a prime of $\bq(\zeta_q^{(3)})$ of local degree $d$ over $\bq$. Then the reduction $A_v$ is an abelian variety over a finite field with Newton polygon \vspace{-0.1cm}$$d\times \frac{j}{d},\; d\times \frac{d-j}{d}$$ for some integer $j$ prime to $d$. Thus, $A_v$ is an abelian variety of type $\IV(1,d)$. The center of $\End^0(A_v)$ is the imaginary quadratic field $\bq(\sqrt{-q})$ and $\End^0(A_v)$ is ramified only at the two places of $\bq(\sqrt{-q})$ lying over $\charac(v)$.\end{Example}

\begin{subsection}{\fontsize{11}{11}\selectfont Arithmetic computations on Jacobians}\end{subsection}

By [Sut18], the cost of addition/doubling of points on a genus three hyperelliptic curve with a Weiertrass point is $I+79M/ 2I+82M$ where $I$ denotes the number of inversions and $M$ the number of multiplications of elements in the finite field of definition. This is the best optimization for point addition/doubling as far as we know. For planar quartic genus three Jacobians, the cost of addition is $163M + 2I$ while that of doubling is $185M + 2I$ ([FOR04]).

\bigskip






\bigskip

\begin{center}\textbf{References} \end{center}
\small

\noindent [AG60] M. Auslander, O. Goldman, \textit{Maximal Orders}, Transactions of the AMS, 1960\\

\noindent [Car03] L. Arenas-Carmona, \textit{Applications of spinor class fields: embeddings of orders and quaternionic lattices}, Ann.Inst. Fourier (Grenoble) 53 (2003), no. 7, 2021–2038\\

\noindent [CCO14] B. Conrad, C. Chai, F. Oort, \textit{Complex Multiplication and Lifting Problems}\\

\noindent [CJS12] Childs, Jao, Soukharev, \textit{Constructing elliptic curve isogenies in quantum subexponential time.}\\

\noindent [CLS11] A. Carocca, H. Lange, R. Rodriguez, \textit{Jacobians with complex multiplication}, Transactions of the AMS, 2011\\

\noindent [Del16] C. Delfs, \textit{Isogenies and Endomorphism Rings of Abelian Varieties of Low Dimension}\\

\noindent [Ell01] J. Ellenberg, \textit{Endomorphism algebras of Jacobians}, Advances in Mathematics 162, 243-271 (2001)\\

\noindent [FOR04] S. Flon, R. Oyono, C. Ritzenhaler, \textit{Fast multiplication of non-hyperelliptic genus $3$ curves}\\

\noindent [FT19] E.V. Flynn, \textit{Genus Two Isogeny Cryptography}, PQCrypto 2019\\

\noindent [Fro68] A. Frohlich, \textit{Formal groups}, Springer-Verlag, Berlin, New York, 1968\\

\noindent [Gal] S. Galbraith, \textit{Supersingular curves in Cryptography}\\

\noindent [Har74] P. Hartung, \textit{Proof of the existence of infinitely many imaginary quadratic fields whose class number is not divisible by 3},  J. Number Theory 6 (1974), 276–278\\

\noindent [Koh96] D. Kohel, \textit{Endomorphism rings of elliptic curves over finite fields}, PhD Thesis\\

\noindent [KLPT], D. Kohel, K. Lauter, C. Petit, J.P. Tignol, \textit{On the quaternion $l$-isogeny problem}\\

\noindent [KMRT98] Knus, Merkujev, Rost, Tignol, \textit{The book of involutions}, American Mathematical Society Col. Publications, 44, AMS 1998\\

\noindent [KO99] W. Kohnen, K. Ono, \textit{Indivisibility of class numbers of imaginary quadratic fields and orders of Tate-Shafarevich groups of elliptic curves with complex multiplication}, Invent. Math. 135 (1999), 387–398.\\

\noindent [La15] K.H.M. Laine, \textit{Security of genus $3$ curves in cryptography}, PhD Thesis\\

\noindent [Mil86] J. S. Milne, \textit{Jacobian varieties},  Arithmetic geometry (Storrs, Conn., 1984), 1986, pp. 167–212\\

\noindent [Mum70] D. Mumford, \textit{Abelian Varieties}\\

\noindent [Oort95] F. Oort, \textit{Abelian Varieties over finite fields} (Notes)\\

\noindent [OU73]  F. Oort, K. Ueno, \textit{Principally polarized abelian varieties of dimension two or three are Jacobian varieties}. J. F. Sci. Univ. Tokyo Sect. IA Math. 20, 377–381 (1973)\\

\noindent [Rei75] I. Reiner, \textit{Maximal Orders}, Academic Press, 1975\\

\noindent [Shim60] G. Shimura, \textit{Abelian Varieties with Complex Multiplication and Modular Functions}, Princeton University Press\\

\noindent [Shio67] T. Shioda, \textit{On the Graded Ring of Invariants of Binary Octavics}, American Journal of Mathematics Vol. 89, No. 4 (Oct., 1967), pp. 1022-1046\\

\noindent [ST66] J.P. Serre, J. Tate, \textit{Good reduction of abelian varieties}, Annals of Math\\

\noindent [Tat69] J. Tate, \textit{Classes d'isogenie des varietes abeliennes sur un corps fini}\\

\noindent [Th17] S. Thakur, \textit{On some abelian varieties of type $IV$}, Preprint\\

\noindent [Th18] \underline{\;\;\;\;\;\;\;}, \textit{Abelian varieties in pairing-based cryptography}\\

\noindent [Voi15] J. Voight, \textit{Quaternion algebras}\\

\noindent [Was82] L. Washington, \textit{Introduction to cyclotomic fields}\\

\noindent [Wat69] W. Waterhouse, \textit{Abelian varieties over finite fields}, Ann. Sci. École Norm. Sup. (4), 2:521–560, 1969\\

\noindent [Yu11] C.F. Yu, \textit{On the existence of maximal orders},  Int. J. Number Theory , 7(2011), 2091-2114\\

\bigskip

\normalsize
\noindent Email: stevethakur01@gmail.com
\bigskip

\noindent Steve Thakur\\
Cybersecurity Research Group\\
Battelle Institute\\
Columbus\\
OH, USA.


\end{document}